\documentclass[12pt,reqno]{amsart}
\usepackage{geometry}                		
\usepackage{graphicx}				
\usepackage{amssymb}
\usepackage{epstopdf}
\usepackage{datetime}
\usepackage{currfile} 

\usepackage{color}

\allowdisplaybreaks

\newtheorem{thm}{Theorem}[section]

\newtheorem{lem}[thm]{Lemma}
\newtheorem{prop}[thm]{Proposition}
\theoremstyle{definition}

\theoremstyle{remark}
\newtheorem{rem}{Remark}[section]

\numberwithin{equation}{section}
						
\newcommand{\R}{\mathbb{R}}
\newcommand{\eps}{\varepsilon}
\newcommand{\del}{\partial}
\newcommand{\cE}{\mathcal{E}}

\renewcommand{\div}{{\rm div}}

\numberwithin{equation}{section}

\newcommand{\red}[1]{\textcolor{black}{#1}}
\newcommand{\blue}[1]{\textcolor{black}{#1}}
\newcommand{\tcc}[1]{\textcolor{black}{#1}}

\newcommand{\newred}[1]{\textcolor{black}{#1}}

\usepackage[nodisplayskipstretch]{setspace}

\usepackage{xpatch}
\makeatletter
\patchcmd{\@tocline}{\hfil}
{\nobreak\leaders\hbox{\ifnum#1<2\hfill\else$\m@th%
\mkern 4.5 mu\hbox{.}\mkern 4.5 mu$\fi}\hfill\nobreak}{}{}
\def\l@section{\@tocline{1}{10pt}{1pc}{}{\bfseries}}
\def\l@subsection{\@tocline{2}{0pt}{\dimexpr 1pc+2em}{}{}}
\makeatother

\begin{document}

\title[]
      {Variational structure of Fokker-Planck equations with variable mobility}
       \author{ Hailiang Liu and Athanasios E. Tzavaras}
\address{Iowa State University, Mathematics Department, Ames, IA 50011} 
\email{hliu@iastate.edu}
\address{Computer, Electrical, Mathematical Sciences and Engineering Division, King Abdullah University of Science and Technology (KAUST), Thuwal, Saudi Arabia.}
 \email{athanasios.tzavaras@kaust.edu.sa}

\date{\today} 


\begin{abstract} 
\tcc{We introduce a generalized Wasserstein metric associated for Fokker-Planck equations with spatially varying mobility.  Defined via  an optimal transport problem,  the cost directly incorporates the frictional dissipation generated by the mobility tensor. By viewing the geometry through a Riemannian lens and applying Nash's embedding theorem, we establish an explicit characterization of optimal transport maps and the underlying variational structure. Under this framework,  
we prove that  the  Fokker–Planck equations arises as a  gradient flow.  Finally, we construct a Jordan–Kinderlehrer–Otto (JKO) time-discrete scheme and establish its convergence to weak solutions.}

 \end{abstract}

\maketitle

\baselineskip 12pt

\section{Introduction}
This work is concerned with Fokker-Planck equations with spatially varying mobility of the form 
\begin{equation} \label{main}
\begin{aligned}
\partial_t \rho &= \nabla_x\cdot \left(\rho A(x) \nabla_x  \delta_\rho F \right),  \quad \text{in } (0,\infty)\times \mathbb{R}^d,
\\
 \rho(0, \cdot ) & = \rho_0( \cdot ) \quad \text{in } \mathbb{R}^d,
 \end{aligned}
\end{equation}
where $A(x)$ is a symmetric positive definite matrix encoding heterogeneous anisotropic diffusion, \red{satisfying for some $c, C > 0$
\begin{equation}\label{hypoA}
c I \le A(x) \le C I \, ,  \quad x \in \R^d \, .
\end{equation}
$F(\rho)$ is an energy functional with $\delta_\rho F = \frac{\delta F}{\delta \rho}$ denoting the first variation of $F$.} The initial $\rho_0$  is a probability density on $\mathbb{R}^d$,  and the solution $\rho(t, \cdot)$ is sought as a time-dependent probability density  satisfying $\rho\geq 0$ almost everywhere in $(0,\infty) \times \mathbb{R}^d$ and $\int_{\mathbb{R}^d} \rho dx = 1$ for $t\geq 0$.

\blue{The primary objective of this work is to establish the optimal transport geometry associated with a variable mobility tensor A(x). Specifically, we introduce a generalized Wasserstein distance induced by the corresponding Riemannian metric,
  providing a direct characterization of its optimal transport maps, variational formulation, metric properties, and dynamic representation. This construction generalizes the classical quadratic-cost Wasserstein distance while maintaining a clear intrinsic interpretation through Riemannian geometry. }
  
\blue{As an application, we show that the Fokker–Planck equation can be framed as a gradient flow with respect to this generalized metric. While gradient-flow formulations for variable mobility have previously been developed within abstract frameworks, most notably by Lisini \cite{Li09},  our focus  is distinct. Rather than concentrating on abstract partial differential equations existence theory, we construct the transport metric explicitly, establish its geometric and variational properties, characterize the associated optimal transport maps, and demonstrate how the gradient-flow structure emerges naturally from this geometric framework.}

Fokker--Planck equations play an important role in modeling fluctuations of physical and biological systems and are
 closely connected with the theory of stochastic differential equations. 
 A classical example of \eqref{main} arises when the energy functional is given by $F=\int \rho {\rm log} \rho dx + \rho \Psi(x) dx$ and $A(x)=I$,   reducing to the usual linear Fokker-Planck equation,
\begin{equation} \label{FP}
\partial_t \rho =\nabla_x\cdot \left( \nabla_x \rho +\rho \nabla_x \Psi \right) \, , 
\end{equation}
which frequently appears in probability theory and is in connection to various optimization problems. Solutions of \eqref{FP} can be constructed via 
the Wasserstein steepest descent method \cite{JKO98, Ot01} (often called JKO-scheme) which proceeds
by solving successively  the  minimization scheme
 \begin{align}\label{JKO}
\rho^{n+1}= \displaystyle{  {\rm Argmin}_{\rho \in K} } \left\{  \frac{1}{2 \tau}W(\rho, \rho^n)^2 +F(\rho)\right\},
 \end{align}
where $K$ denotes the set of probability densities on $\mathbb{R}^d$ with finite second moments and the iteration is initialized with $\rho^0=\rho_0$. \red{Here, $\tau >0$ is the time step.}  The variational construction interprets  (\ref{FP}) as steepest descent of the functional $F(\rho)$ with respect 
to the 2-Wasserstein metric. The approach has been applied to a wide class of energy functionals and 
extended to more complex systems with analogous gradient flow structures, see, for example,  
\cite{JKO98,Ot01,BCC08,BL13,MMS09,Sa15,PMSV16,KMX17}.
 
 The studies  \cite{JKO98, Ot01}  attracted considerable attention because they introduced a new framework for advection-diffusion
 equations but also due to the connection with gradient flows on Riemannian manifolds. 
 Motivated by the De Giorgi \cite{Gi93} minimizing movement scheme,
  Ambrosio-Gigli-Sava\'{r}e \cite{AGS08} developed a theory of curves of maximal slope in metric spaces of probability
  measures leading to a gradient flow, 
\begin{equation}\label{riemflow}
 \dot x(t) = - {\rm grad}_d F(x(t)) \, ,
\end{equation}
on a Riemannian manifold of probability measures. 
For various aspects of that theory,  the reader is referred to \cite{Ot01,Br03,Vi03,AGS08,OW06,Sa15}.  
Lisini \cite{Li07, Li09} gave a construction of solutions to \eqref{main} using the framework of \cite{AGS08}, introducing
the Wasserstein distance defined in \eqref{rv+} below.

The objective of this work is to complement the approach of Lisini  \cite{Li09} by emphasizing 
the variational structure underlying the heterogeneous diffusion \eqref{main}.  
The central object of interest is the metric 
\begin{equation}\label{rv+}
W_A(\rho_0, \rho_1)^2:= \min_{(\rho, v) \in \mathcal{S}} \left\{  \tau  \int_0^\tau \int_{\mathbb{R}^d} \rho(s, x)v(s, x)\cdot A^{-1} (x) v(s, x)dx ds\right\}, 
\end{equation}
where the minimization occurs over $\mathcal{S}$ the set of all admissible flows $(\rho, v)$  satisfying 
\begin{align*}
\mathcal{S} = \left \{ (\rho,v) :(0,\tau) \times \R^d \to \R^+ \times \R^d  \quad \big |  \quad  
\begin{cases} \quad \partial_s \rho +\nabla_x \cdot(\rho v)=0, & \\
             \rho(0, \cdot)=\rho_0, \quad \rho(\tau, \cdot)=\rho_1  & \\
              \end{cases}\right \} \, .
\end{align*}
\red{Extending the Benamou-Brenier formula  \cite{BB00}, the action functional  \eqref{rv+} reflects  the frictional dissipation of transport through a spatially varying medium, with  
dissipation rate 
$$
\frac{d}{dt} F(\rho)=-\int_{\R^d} \rho v^\top A^{-1}(x)v dx, \quad v = - A \nabla_x \delta_\rho F.
$$
As $\tau>0$, the JKO minimizing movement scheme
\[
\rho^{n+1}={\rm Argmin}_{\rho \in \mathcal K}  \left \{ \frac{1}{2\tau} W_A^2(\rho, \rho^n) + F (\rho) \right\}
\]
yields the discrete velocity 
\[
\frac{r^{n+1}(x) - x}{\tau}
\approx
- A \nabla_x  \delta_\rho  F (\rho^{n+1}),
\]
establishing  (\ref{main}) as a $W_A$-gradient flow.  This reflects a maximal energy dissipation principle \cite{AGS08,Li09} and provides a constructive Eulerian framework for optimal maps  \cite{BB00, Sa15}.}

The analysis of the minimization problem \eqref{rv+} is based on viewing the friction matrix $B(x) : = A^{-1}(x)  >0$ as a Riemannian metric.
The Nash-Kuiper theorem \cite{Na56, GR70} guarantees the decomposition of the form 
\begin{align}\label{abb-}
B (x)=(\nabla_x b)^\top (\nabla_x b)
\end{align}
for a smooth injective map $b: \mathbb{R}^d \to \mathbb{R}^q$ with $b(0)=0$, $q> d$.  Using \eqref{abb-},  we show that the above 
variational formulation is equivalent to the classical  Monge problem, 
\begin{align}\label{rss+}
r^*(x)={\rm arginf}_{r\in \mathcal{A}} \int d_M^2 \big( b(x),  b(r(x) \big) \, \rho_0(x)dx,
\end{align}
over the class  $ \mathcal{A}$ of admissible measure-preserving maps pushing $\rho_0$ to $\rho_1$, where  
$d_M$ denotes the intrinsic Riemannian distance on $M=b(\mathbb{R}^d)$.  


For the classical Wasserstein distance,  the optimal transport theory is well developed for translation-invariant costs of the form $h(x-y)$ with strictly convex  $h$ \cite{Br91, GM95, Vi03, Sa15}. 
\blue{The framework extends that perspective to transport costs governed by non-Euclidean geometry. By  lifting the problem through 
an appropriate Riemannian embedding,  we analyze the intrinsic cost through an optimal transport problem in the ambient space. This connection allows us to transfer core results from classical theory
\cite{Br91,Ev89, GM96, RR98, Mc2001} to the variable-mobility setting. } Related results on convexity and contraction properties in Wasserstein spaces can be found in \cite{M97, OW06}. 

In the final part, we construct a time discrete, iterative variational scheme, 
 \begin{align*}
\rho^{n+1}= \underset{\rho \in K}{\rm argmin} \left\{  \frac{1}{2 \tau}W_A(\rho, \rho^n)^2 +F(\rho)\right\} \, ,
  \end{align*}  
initiating with $\rho^0 = \rho_0$. This minimizing scheme leads to a constructive existence proof for weak solutions to (\ref{main}) which can be 
numerically implemented
and confirms that (\ref{main}) is a gradient flow for the free energy functional $F(\rho)$ in the weighted Wasserstein metric $W_A$. \red{While the variational time discretization admits a general formulation, the  convergence  results obtained here are restricted to a particular energy functional $F$}.   One ingredient is to carefully trace how the mobility matrix $A$ couples the Riemannian point of view and the 
classical partial differential equation framework. The variational scheme yields an approximation of the implicit Euler scheme
$$
\frac{\rho^{n+1}(x)-\rho^n(x)}{\tau}=\nabla_x \cdot ( \rho^{n+1} A(x) \nabla_x \delta_\rho F|_{\rho=\rho^{n+1}}), 
$$
with the error controlled quadratically in the weighted Wasserstein distance.  We show that, under some
strong compactness assumptions, $\rho$ is a weak solution to (\ref{main}).  

To focus on the transport geometry, we analyze the linear Fokker–Planck equation and defer nonlinear compactness challenges to future work (cf. \cite{BL13,KMX17,Li09,PMSV16,Sa15, FT22, EHS25}). 
\blue{While closely related to the generalized JKO framework for convex Hamiltonians \cite{FGY11}, our mobility-induced Hamiltonian in (\ref{main}) yields explicit transport metrics, map characterizations, and variational gradient-flow structures, thus complementing \cite{FGY11}.}

The paper is organized as follows. Section \ref{sec:gradient}  introduces the class of Fokker--Planck equations with heterogeneous anisotropic diffusion,  
discusses their gradient-flow structure, and explains their relation to the  small-mass limit  of damped Euler systems.  In Section \ref{sec:wasserstein},
we develop the generalized Wasserstein metric, establish its main properties and equivalent formulations,  and prove the existence and uniqueness of the associated optimal transport maps.
Section \ref{sec:variation} develops the formal variational limit and identifies the limiting equation, while Section \ref{sec5} provides the existence theory for the minimizing scheme and completes the convergence proof in the case of linear diffusion. 


In Appendix \ref{sec:appA},  we derive the de-Giorgi minimizing movement scheme as a high-friction limit 
of variational approximations of second-order evolutions in Lagrangian coordinates. 
Appendix \ref{sec:isoemb}  collects the background material on the isometric embedding results used in the paper. 

\section{Gradient flow and the high-friction limit}\label{sec:gradient}

For mechanical systems with friction in Lagrangian coordinates,  it is part of the folklore of the subject
how the De Giorgi minimizing movement method emerges in the small mass approximation
of variational approximations as a relaxation problem; we refer to the appendix for a presentation of its ramifications.
We are here interested in how the  JKO scheme of  Jordan-Kinderlehrer-Otto \cite{JKO98}, Otto \cite{Ot01}
and the fluid mechanics interpretation of Wasserstein of Brenier-Benamou \cite{BB00} have to be 
adapted for the case of heterogeneous-anisotropic diffusion. \red{ Most of the material in this section is standard and is presented only to provide background and fix notation. 
One notable exception is the remarks in Section \ref{sec:maxenedis}  which illustrate the role of frictional dissipation in Euler systems
to the form of the Wasserstein distance for Fokker-Planck equations.
The derivations  in this section are formal requiring smooth solutions.} 
\subsection{Gradient flow structure} \label{Sec2.1}
In this section, we review the geometric framework of \cite{Ot01}
and interpret \eqref{main} as a gradient flow.  
Equation  \eqref{main} can be written as a continuity equation 
\begin{equation} \label{main+}
\partial_t \rho +\nabla_x\cdot \left(\rho u\right)=0 \, ,
\end{equation}
for the mass density  $\rho$ and the velocity $u$, coupled with Darcy's law $u = -A(x) \nabla_x p$. 
Here, $p$ is the pressure and $A(x)$ is a symmetric positive definite matrix
describing the mobility of a non-homogeneous and anisotropic medium.  
The  pressure is determined thermodynamically through the chemical potential
 \begin{equation} \label{pe}
p=\delta_\rho F,
\end{equation}
where $F$ is the free energy functional.  The geometric formulation introduced in  \cite{Ot01} considers the space of probability densities 
$$
\mathcal{M}=\left\{ \rho \geq  0, \quad \int_{\mathbb{R}^d}  \!\!\rho dx =1\right\} \, ,
$$
equipped with a Riemannian metric $g$ and an energy functional  $F(\rho)$.  
The metric tensor converts the differential $\delta_\rho F$
which is a cotangent vector field onto the gradient $\rm{grad} F$ which is a tangent field, so that 
 for all vector fields $s$ on $ \mathcal{M}$ we have at the point $\rho \in \mathcal{M}$:
\begin{equation}\label{2.6}
g_\rho ( \rm{grad} F , s ) = \langle \delta_\rho F , s \rangle = \int_{\mathbb{R}^d}  \delta_\rho F \,  s  \, dx.
\end{equation}
The gradient flow 
\begin{equation}
\del_t \rho = - \rm{grad} F \big |_\rho
\end{equation}
is interpreted as 
$$
g_\rho (\del_t \rho , s) + \int_{\mathbb{R}^d}   \delta_\rho F  s dx = 0. 
$$
This formulation immediately yields the energy dissipation identity
$$
\frac{d}{dt}F(\rho)=\langle {\delta_\rho F}, \partial_t \rho\rangle =-g_\rho(\partial_t \rho, \partial_t \rho). 
$$
For the Fokker-Planck equation (\ref{main}), the evolution preserves both the non-negativity of $\rho$ 
and the total mass $\int \rho dx $. 

The tangent space of  $\mathcal{M}$ is given by 
$$
T_\rho \mathcal{M}=\left\{ s \;  \;  \int_{\mathbb{R}^d} \!\! s dx =0\right\}.
  $$
A tangent vector can be represented through a potential  $p$ on  $\mathbb{R}^d$, defined up to an additive constant,  by the elliptic relation 
 \begin{equation} \label{sp}
  s=-\nabla_x \cdot (\rho A \nabla p).
  \end{equation}
  For bounded domains, the zero-flux boundary condition is imposed.

The Riemann metric  is then defined by
 \begin{equation}\label{defmetric}
 \begin{aligned}
 g_\rho(s_1, s_2) &=\int \rho \nabla_x p_1 \cdot A \nabla_x p_2 dx,
 \\
\mbox{where $p_i$ and $s_i$ are related via} \quad &s_i = -\nabla_x \cdot (\rho A \nabla p_i ) \, , \quad i = 1,2. 
\end{aligned}
\end{equation}
Note that  $g_\rho (s_1, s_2) = \int_{\mathbb{R}^d}  s_1 p_2 dx$.
Using (\ref{2.6}) with $s=\partial_t \rho$, we obtain    
$$
\int_{\mathbb{R}^d}  (\partial_t \rho) \,  p dx +\int_{\mathbb{R}^d}   {\delta_\rho F} s dx =0,
$$
where $p$ is related to $s$ via $-\nabla_x \cdot  (\rho A\nabla  p)=s$. After an integration by parts this gives 
$$
\int_{\mathbb{R}^d}  \big (\partial_t  \rho - \nabla \cdot (\rho A\nabla  \delta_\rho  F ) \big ) p dx=0,
$$
which recovers the Fokker-Planck equation (\ref{main}). 

\medskip
\noindent
{\bf Properties of the metric}.
Let  $B(x) := A^{-1}(x)$ which represents the friction of the medium. The metric admits the following equivalent formulations.
\begin{itemize}
\item[(i)]
 If  $v_i = A \nabla_x p_i$,  then $s_i = - \div (\rho v_i)$ and
$$
g_p (s_1 , s_2 ) = \int_{\mathbb{R}^d}  \rho v_1 \cdot B v_2 \, dx. 
$$

\item[(ii)]  \red{For fixed $ \rho\in\mathcal M$ and $ s\in T_\rho\mathcal M$},  the metric can be
characterized by the minimization problem
\begin{align}
& \min_v
\left \{ 
\int_{\mathbb{R}^d}  \rho v\cdot B  v  dx  \; \Big| \; \text{subject to }
 \; \;  s + \nabla_x  \cdot(\rho v)=0  \right  \}.
\label{2.8} 
\end{align} 
The unique minimizer is  $v = \rho A \nabla p$, where  $p$ solves  \eqref{sp}. 
\end{itemize}

Indeed, introducing $m = \rho v$, the problem becomes   
$$
\min_{m } \int_{\mathbb{R}^d}  \frac{1}{\rho} m \cdot B m dx
$$
subject to $ - \nabla_x \cdot  m= s$.
The Euler Lagrange equations for the augmented Langrangian gives 
$$
\frac{1}{\rho} B m  - \nabla_x \varphi = 0,
$$
and hence  $m = \rho A \nabla_x \varphi$. The constraint then yields  
 $$
 - \nabla_x \cdot  ( \rho A \nabla_x \varphi) = s,
 $$
 which proves the characterization.

\subsection{Fokker-Planck with variable mobility as a high-friction limit}\label{sec:FPvar}
Now (\ref{main}) is viewed as a high-friction limit of an Euler system.
Consider the (rescaled in time) hydrodynamic system 
\begin{subequations}\label{1.1}
\begin{align}
& \partial_t \rho +\nabla_x \cdot (\rho u)=0,\\
& \eps \, \rho D_t u =  -  \rho \nabla_x h^\prime (\rho)  - \rho B (x)u \, ,
\end{align}
\end{subequations}
where $D_t: = \partial_t  +u\cdot \nabla_x$ is the material derivative,  $\eps >0$ is a small mass parameter, 
 $B := A^{-1}$ is a symmetric,  positive definite matrix,  and consider a potential energy functional
 $$
 F[ \rho ] = \int h(\rho) dx \, .
 $$
 The $L^2$ first variation of $F[\rho]$ is computed as $\delta_\rho F = h^\prime (\rho)$. The resulting system is an Euler system with friction
 described by $B(x)$. It satisfies the energy dissipation identity
 \begin{equation}\label{enedis}
\frac{d}{dt}  \left( F(\rho)+ \int_{\mathbb{R}^d}   \frac{\eps}{2} \rho|u|^2dx  \right)+\int_{\mathbb{R}^d}  \rho B (x)u\cdot u dx =0, 
\end{equation}
In the limit $\epsilon \to 0$, \eqref{1.1} formally yields \eqref{main},
\begin{equation}\label{disprob}
\begin{aligned}
 \partial_t \rho +\nabla_x\cdot (\rho u) &=0, 
 \\
 u&= - A(x) \nabla_x \delta_\rho F \, ,
\end{aligned}
\end{equation}
which inherits  from \eqref{enedis} the potential energy dissipation identity
\begin{align}\label{potdis}
\frac{d}{dt} F(\rho) + \int_{\mathbb{R}^d}   \rho \nabla_x  \delta_\rho F \cdot A(x) \nabla_x \delta_\rho F  dx  = 0 \, .
\end{align}
\subsection{Maximal energy dissipation}\label{sec:maxenedis} We record an interesting connection between the problems \eqref{1.1} and \eqref{disprob}.
Consider a flow $(\rho, u)$ that satisfies the mass balance \eqref{main+}  occuring in a frictional landscape characterized
by the friction matrix $B(x)$. Consider  the  balance of total energy
$$
\cE (t) = \int_{\mathbb{R}^d }\frac{\eps}{2} \rho |u|^2 \, dx + F(\rho), 
$$
where $\eps$ is a mass parameter, and $F(\rho)$ a functional determining the potential energy.  Using \eqref{main+} and the association $B = A^{-1}$, we compute
$$
\begin{aligned}
\frac{d}{dt} \cE (t)  &=\int_{\mathbb{R}^d }  \eps  ( \rho_t  |u|^2/2  + \rho u u_t ) \, dx + \langle \delta_\rho F , \rho_t \rangle
\\
&= \int_{\mathbb{R}^d  } \rho \left ( \eps D_t u + \nabla_x \delta_\rho F \right ) \cdot u \, dx
\\
&= \int_{\mathbb{R}^d } \rho  \Big \langle \eps D_t + \nabla_x \delta_\rho F  \; ,  \; B(x)   u \Big \rangle_{A(x)}  \, dx
\\
&= \frac{1}{2}  \int_{\mathbb{R}^d  } \rho \Big | \eps D_t + \nabla_x \delta_\rho F + B(x) u \Big |^2_{A(x)} dx
\\
&\qquad - \frac{1}{2}  \int_{\mathbb{R}^d  } \rho \Big | \eps D_t + \nabla_x \delta_\rho F  \Big |^2_{A(x)} dx 
- \frac{1}{2} \int_{\mathbb{R}^d } \rho u \cdot B(x) u \, dx, 
\end{aligned}
$$
where we used the notation $\langle a, b \rangle_{A(x)} = a \cdot A(x) b$ for the inner product associated to the matrix $A(x)$ 
(which is the natural product when $A(x)$ is viewed as a Riemannian metric),
$|a|_A^2 =  a\cdot Aa$ the corresponding norm, and the formula $|B(x) u|^2_{A(x)} = u \cdot B(x) u$.

We conclude that along a mass conserving flow $(\rho, u)$ on a frictional landscape with friction $B(x)$ the rate of total energy decay
has a lower bound
$$
\frac{d}{dt}  \Big (  \int_{\mathbb{R}^d  }\frac{\eps}{2} \rho |u|^2 \, dx + F(\rho) \Big ) 
\ge 
- \frac{1}{2}  \int_{\mathbb{R}^d  } \rho \Big | \eps \frac{Du}{Dt} + \nabla_x \delta_\rho F  \Big |^2_{A(x)} dx 
- \frac{1}{2} \int_{\mathbb{R}^d } \rho u \cdot B(x) u \, dx. 
$$
The maximal dissipation occurs for the flow $(\rho, u)$ defined by \eqref{1.1} and  satisfies the energy dissipation \eqref{enedis}.
This maximal rate of energy dissipation for the mechanical energy bears a striking resemblance to the maximal (potential) energy dissipation 
for gradient flows, which  plays a significant role in the construction of solutions to gradient flows in Riemannian manifolds,
see \cite{AGS08}  and \cite{Li09}.  The argument presented here further connects the limiting Fokker-Planck equation with the
way it emerges as a small-mass approximation of an Euler flow with friction.
This argument should be compared to the entropy rate admissibility criterion of Dafermos \cite{Daf12} and offers support to that criterion in a special context.
\red{The same computation applies on a bounded domain $\Omega$ with no-flux boundary conditions $u \cdot \nu = 0$ on
$\del \Omega$, or under periodic boundary conditions on the torus.}
\section{A weighted Wasserstein metric and optimal maps} \label{sec:wasserstein}
\subsection{Preliminaries on optimal transportation} 
Let $\mu_0, \mu_1$ be two Borel probability measures on $\mathbb{R}^d$.
Define the admissible class 
$$
\mathcal A: \{ r: \mathbb{R}^d  \to \mathbb{R}^d:\;  r_{\rm \#} \mu_0=\mu_1\},
$$
\red{where $r_{\rm \#}$  denotes the push-forward of measure}.  
The Monge problem seeks an optimal transport 
$$
r^*={\rm arg} \inf_{r\in \mathcal{A}}  \int_{\mathbb{R}^d} c(x, r(x))d\mu_0(x), 
$$
where $c(x, y)$ is the transportation cost.  
In general, the infimum need not be attained. 
Kantorovich relaxed this problem by replacing transport maps with transport plans.  Let 
\begin{equation}\
\label{defpi}
\Pi(\mu_0,\mu_1)
=
\left\{
\pi\in\mathcal P(\mathbb{R}^d\times\mathbb{R}^d):
(\mathrm{proj}_1)_{\#}\pi=\mu_0,\;
(\mathrm{proj}_2)_{\#}\pi=\mu_1
\right\}
\end{equation}
denote the set of probability measures with marginals $\mu_0$ and $\mu_1$. The corresponding Kantorovich problem is
\begin{align}\label{da}
\min_{\pi \in \Pi(\mu_0, \mu_1)}
\int_{\mathbb{R}^d \times \mathbb{R}^d} c(x, y)
 d\pi(x, y). 
\end{align}
This relaxed formulation is considerably more tractable; see \cite{Ev89}. 
Its minimizers are called \emph{optimal transport plans}. A natural question is whether the relaxation changes the optimal value, namely, 
whether the minimum value of the Kantorovich problem coincides with the infimum of the Monge problem. 

If the cost $c(x, y)$ is continuous, and the source measure $\mu_0$ is atomless,  then it is known  \cite{Sa15} that 
$$
\inf_{r\in \mathcal{A}} \int_{\mathbb{R}^d} c(x, r(x))d\mu_0(x)=\min_{\pi \in \Pi(\mu_0, \mu_1)}
\int_{\mathbb{R}^d \times \mathbb{R}^d} c(x, y) d\pi(x, y).
$$
In this setting, an optimizer of the Monge problem  is  called  an optimal transport map,  and the minimizer of the Kantorovich problem is called an optimal transport plan. When the optimal plan is induced by a map $r$, 
it takes the form $\pi=(id, r)_{\#}\mu_0$.

For the quadratic cost $c(x, y)=|x-y|^2$,  let $d\mu_0 = \rho_0 dx$, $d\mu_1 = \rho_1 dx$, where  $\rho_0$, $\rho_1$ are smooth probability densities.
\red{A fundamental result (see \cite{Br91,Ev89, GM96, RR98})  asserts that there is a unique optimal transport $r$ characterized as the
unique map pushing$\rho_0$ to $\rho$.}  Moreover, this map is characterized by the fact that it can be written as the gradient of a convex function,
$$
r(x)=\nabla_x u(x). 
$$
The potential $u$ satisfies, in a suitable weak sense, the Monge-Amp\'{e}re equation 
$$
{\rm det} (\nabla_x^2 u(x))\rho_1(\nabla_x u(x))=\rho_0(x). 
$$
Under additional geometrical assumptions on the domain and the densities, the smoothness of $\rho_0$ and $\rho$ is inherited by the potential $u$, as established by Caffarelli
\cite{Ca96+}.

\subsection{The analog of Wasserstein for heterogeneous, anisotropic mobility}
\red{ The objective is to identify a Riemannian distance that plays the role of  Wasserstein in a
heterogeneous and anisotropic landscape. Recall that $A(x)$ is a smooth positive definite and symmetric matrix-valued function
and  in view of \eqref{hypoA} its inverse $B(x) = A^{-1} (x) $ satisfies  for some $c_0 > 0$}
\begin{align}\label{c0}
B(x) \ge c_0 I \, .
\end{align}
We posit that $B $ admits the decomposition 
\begin{align}\label{abb}
B (x)=(\nabla_x b)^\top (\nabla_x b)
\end{align}
for a smooth injective map $b: \mathbb{R}^d \to \mathbb{R}^q$ with $q\geq d$.  
\tcc{ We place here \eqref{abb} as an assumption,
but note that such a relation is ensured in various contexts by the Nash embedding theorem \cite{Na56, GR70},
which provides an isometric embedding into $\mathbb{R}^q$, for some large $q$.}  We refer to  Appendix \ref{sec:isoemb} for details on the isometric embedding literature. Various  aspects of
the forthcoming analysis exploit this geometric property.

Inspired by \red{ Section \ref{Sec2.1}} and \cite{BB00},  we propose the following definition for  $W_A (\rho_0, \rho_1)$: 
Fix an interval $[0,  \tau]$ and consider all sufficiently smooth density and velocity fields $\rho(s, x) \geq 0$, $v(s, x)\in \mathbb{R}^d$
satisfying the target problem
\begin{align}\label{tr}
\partial_s \rho +\nabla_x \cdot(\rho v) &=0 \, ,  \qquad (s,x) \in (0,\tau) \times \R^d \, ,
\\
\label{r01}
\rho(0, \cdot)=\rho_0, \quad \rho(\tau, \cdot)&=\rho_1, 
\end{align}
subject to the initial and final conditions \eqref{r01}.  
Let $\mathcal{S}$ denote the class of admissible fields $(\rho, v)$,
and define
\begin{equation}\label{rv}
W^2_A(\rho_0, \rho_1): =  \min_{(\rho, v) \in \mathcal{S}} \tau \int_0^\tau \int_{\mathbb{R}^d} \rho(s, x)v(s, x)\cdot B (x) v(s, x)dx ds. 
\end{equation}
 Let $b: \mathbb R^d \to \mathbb R^q$ be a smooth embedding, and set $M=b(\mathbb{R}^d)$. Then $M$  is a smooth d-dimensional embedded submanifold of $\mathbb R^q$, and $b: \mathbb R^d \to M$ is a diffeomorphism.
 We equip  $M$ with the Riemannian metric induced by the embedding, whose intrinsic distance 
\begin{equation}\label{dm}
d_M(\xi_0, \xi_1)= \left( \inf_{\xi \in AC([0, \tau];  M)}  \tau \int_0^\tau |\dot \xi(s)|^2 ds )\right)^{1/2},
\end{equation} 
subject to the boundary conditions  $\xi(0)=\xi_0$ and  $\xi(\tau)=\xi_1$, where $AC([0, \tau]; M)$ denotes the set of absolutely continuous curves in $M$.  By a standard change of variables, this definition is independent of $\tau>0$, 
and coincides with the usual action formulation parameterized over $[0,1]$.

Consider the Monge–Kantorovich optimal transportation problem 
\begin{equation}\label{opttrans}
 \min_{\pi \in \Pi(\mu_0, \mu_1)}
\int_{\mathbb{R}^d \times \mathbb{R}^d} d_M(b(x), b(y))^2 d \pi(x, y) \, ,
\end{equation}
where $\mu_0 = \rho_o dx$, $\mu_1 = \rho_1 dx$, 
with quadratic cost on the geodesic distance.
This problem has been extensively analyzed in the literature 
(see,  McCann \cite{Mc2001} for compact manifolds, and  Figalli \cite{Fig07},  Figalli and  Fathi \cite{FF10} for non-compact manifolds), and, under an appropriate absolute continuity assumption on the source measure,
has a unique optimal map that admits the Riemannian representation
\begin{equation}
T(\xi)=\exp_\xi(-\nabla_M\phi(\xi)) \, .
\end{equation}
Here $\phi$ is a $c$-convex Kantorovich potential and $\exp$ denotes the Riemannian exponential map; see \cite{Mc2001,FF10}. Pulling back gives 
$$
r^*(x)=b^{-1}\circ T\circ  b(x)=b^{-1}\circ \exp_{b(x)} (-\nabla_M \phi(b(x))).
$$
This gives a geometric characterization of the transport problem.

In addition, this analysis allows to connect the formula \eqref{rv} through the Nash-Kuiper isometric embedding 
Theorem with the quadratic cost problem \eqref{opttrans} through a formal argument:

\begin{prop}\label{pr1}
Under hypotheses \eqref{c0}, \eqref{abb}, 
the minimization problem \eqref{rv} defines
a generalized Wasserstein distance  that satisfies
\begin{equation}
\label{defhetwas+}
W^2_A(\rho_0, \rho_1) = \min_{\pi \in \Pi(\mu_0, \mu_1)}
\int_{\mathbb{R}^d \times \mathbb{R}^d} d_M (b(x), b(y))^2 d \pi(x, y), 
\end{equation} 
with $d\mu_0 = \rho_0 dx$, $d\mu_1 = \rho_1 dx$ and $\Pi(\mu_0, \mu_1)$ defined in \eqref{defpi}. 
\end{prop}

\begin{proof} 
Following \cite{BB00}, we work in
Lagrangian coordinates. Let $\rho_0,\rho_1$ be compactly supported probability densities on $\mathbb R^d$, and let
$(\rho,v)$ be a sufficiently smooth solution of \eqref{tr}--\eqref{r01}. Define the
characteristic flow $X(s,a)$ by
\begin{equation}\label{xv}
X(0,a)=a,\qquad
\partial_sX(s,a)=v(s,X(s,a)).
\end{equation}
Then the continuity equation implies the pushforward relation
\[
\rho(s,\cdot)=X(s,\cdot)_\#\rho_0,
\]
or equivalently, for every test function $f$,
\[
\int_{\mathbb R^d}f(x)\rho(s,x)\,dx
=
\int_{\mathbb R^d}
f(X(s,a))\rho_0(a)\,da.
\]
In particular,
\[
r(a):=X(\tau,a)
\]
is a measure-preserving map from $\rho_0$ to $\rho_1$. Using the above change of variables together with
\eqref{abb},
\[
B=(\nabla b)^\top\nabla b,
\]
we obtain
\begin{align*}
&\tau
\int_0^\tau\!\!\int_{\mathbb R^d}
\rho\,v^\top Bv\,dx\,ds   =
\tau
\int_0^\tau\!\!\int_{\mathbb R^d}
\rho_0(a)\,
\bigl|\partial_sb(X(s,a))\bigr|^2
\,da\,ds.
\end{align*}
By the definition of the Riemannian distance, this is bounded from below by
$$
\int_{\mathbb R^d}
d_M^2(b(r(a)), b(a))
\rho_0(a)\,da.
$$
The reverse inequality follows immediately by choosing the flow generated by an 
optimal Monge transport map $r=r^\star$, for which each trajectory is a minimizing geodesic connecting $b(a)$ and 
$b(r^\star(a)))$. 
The existence of such an optimal transport map is outlined in Theorem  \ref{thm3.3-riemannian}. 
\end{proof}

Another insight to the quantity
$
\frac{1}{2\tau} W_A(\rho_0, \rho_1)^2
$
can be obtained by the optimal control problem:
\begin{align}\label{op} 
& {\rm inf}_{\rho,  v} \left( \frac{1}{2} \int_0^\tau \int_{\mathbb{R}^d} \rho(s, x)v(s, x)\cdot B (x) v(s, x)dxds \right)\\ \notag
& \text{subject to } \quad \partial_s \rho +\nabla_x \cdot(\rho v)=0,\\ \notag
& \qquad \qquad \rho(0, \cdot)=\rho_0, \quad \rho(\tau, \cdot)=\rho_1. \notag 
\end{align}
The Eulerian formulation (\ref{op}) of the optimal control problem in density space  
can be written as a saddle-point problem by introducing
a space-time dependent Lagrange multiplier $\phi(s, x)$ \cite{BB00}. The Lagrangian with $m=\rho v$  is given by 
\begin{align}\notag
\mathcal{L}(\phi, \rho, m)& =\int_0^\tau \int_{\mathbb{R}^d} \left( \frac{1}{2\rho} m\cdot B (x) m -\partial_s \phi \rho -\nabla_x \phi \cdot m\right)dxds \\
 & \qquad  -\int_{\mathbb{R}^d} (\phi(0, x)\rho_0 (x) -\phi(\tau, x)\rho_1(x))dx,
\end{align}
where the terms involving $\phi$ come from the integration by part of the transport equation (\ref{tr}) using boundary conditions (\ref{r01}). 
Given initial and final densities $\rho_0$ and $\rho_1$, the optimal control problem (\ref{op}) is equivalent to
the saddle-point problem:
\begin{align}
{\rm inf}_{\rho, m} {\rm sup}_{\phi} \mathcal{L} (\phi, \rho, m).
\end{align}
 The optimality condition requires that 
  $\delta_m\mathcal{L}=0$ and $ \delta_\rho \mathcal{L} \leq 0$,  which leads to the following property: 
 $$
 m= \rho A(x) \nabla_x \phi, \quad -\frac{1}{2\rho^2}m \cdot B m -\partial_s \phi \leq0.
 $$
 This says that $\phi$ is a super-solution to the Hamilton-Jacobi equation. Here if $\rho>0$, we obtain $ \delta_\rho \mathcal{L} =0 $, which gives 
 \begin{align*}
& \partial_s \phi + \frac{1}{2 \rho^2} m \cdot B (x) m=0, \quad m= \rho A(x) \nabla_x \phi.
\end{align*}
Notice that $m$  can be eliminated, leading to the Hamilton–Jacobi–continuity system
\begin{subequations}\label{rp}
\begin{align}
& \partial_s \phi + \frac{1}{2} \nabla_x \phi \cdot A(x) \nabla_x \phi=0,\\ 
& \partial_s \rho +\nabla_x \cdot(\rho A(x) \nabla_x \phi)=0.
\end{align}
\end{subequations}
The corresponding optimal velocity field is
$$
v(s, x)=A(x) \nabla_x \phi(s, x). 
$$
The characteristic flow $X(s,x)$ is generated by integrating this velocity field:
$$
\partial_s X(s, x)= A(X(s, x)) \nabla_x \phi(s, X(s, x)). 
$$

 \red{
\begin{rem}\label{remcomparison}
(Relation to \cite{Li09}). While both our formulation and that of Lisini \cite{Li09} build upon the friction-matrix metric $B=A^{-1}$,  Lisini works within an abstract Riemannian setting. 
By applying Nash's embedding theorem, we realize this metric extrinsically via an embedding $b: \mathbb R^d \to \mathbb R^q$.  
 This geometric perspective unlocks an explicit gradient characterization for the optimal transport map and yields stronger structural control over the associated Wasserstein distance.
\end{rem} 
}

\subsection{Existence of optimal maps}\label{sec:exiopt}
The existence of optimal transport maps for quadratic cost in the geodesic distance is discussed in 
\cite{Mc2001, Fig07, FF10}. The question is whether an optimal transport plan 
$\pi$ can be represented as induced by a map, namely,  $\pi=(\mathrm{id},  r)_{\#} \mu_0$
for some transport map $r$. 
The dual functional is introduced
$$
J[\phi, \psi]:=\int_{\mathbb{R}^d} \phi(x) \rho_0(x)dx +\int_{\mathbb{R}^d} \psi(y) \rho_1(y)dy,  
$$ 
and the dual problem is solved:  find an optimal pair $(\phi^*, \psi^*)\in \mathcal{L}$ such that 
\begin{align}\label{kan}
J[\phi^*, \psi^*]=\max_{(\phi, \psi)\in \mathcal{L}} J[\phi, \psi],
\end{align}
where 
$$
\mathcal{L}=\{(\phi, \psi):  \;  \mathbb{R}^d \to \mathbb{R} \; \text{continuous}, \; \phi(x) +\psi(y) \leq  c(x, y), x, y \in \mathbb{R}^d\}.
$$
The cost functional and its dual functional $J$ are connected via 
\begin{align}\label{WK}
W_A(\rho_0, \rho_1)^2 =\max_{(\phi, \psi)\in \mathcal{L}} J[\phi, \psi],
\end{align}
which can be ensured to hold if $c$ is continuous (lower semi-continuity is sufficient) and  bounded from below; see \cite[Theorem 1.42]{Sa15}.

Moreover, the primal Kantorovich problem and its dual formulation  are connected through the following characterization.
\begin{lem} Let $X={\rm supp}(\rho_0)$ and  $Y={\rm supp}(\rho_1)$.
 If $\pi$ is an optimal transport plan and $(\phi^*, \psi^*)$ is an optimal dual pair satisfying the c-transform relations  
 \begin{align*}
& \phi^*(x)=\min_{y\in Y} (c(x, y)-\psi^*(y)),\\
& \psi^*(y)=\min_{x\in X} (c(x, y)-\phi^*(x)),
\end{align*}
then the support of $\pi$ is contained in the set
$$
\{(x,y)\in X \times Y \; \text{such that} \; \phi^*(x) +\psi^*(y)= c(x,y)\}\subset X \times Y
 $$
This set is the $c$-subdifferential associated with the potential potential and provides the geometric structure needed to determine when an optimal plan is induced by a transport map.
\end{lem}

The key challenge is to determine when  the optimal transport relation is single-valued, namely, when the optimal plan is induced by a transport map.  
This is addressed in Brenier \cite{Br91} for a quadratic cost in $\R^n$; in the general Riemannian setting, 
the optimal transport map admits the characterization,
see \cite{Mc2001}, \cite{Fig07},  \cite{FF10} :

\begin{thm}\label{thm3.3-riemannian}
Let \(\mu_0\) and \(\mu_1\) be probability measures on \(\mathbb{R}^d\), and let
$
b:\mathbb{R}^d \rightarrow M\subset \mathbb{R}^q
$
be a smooth embedding into a smooth Riemannian manifold \(M\), equipped with the induced Riemannian metric. Define the transport cost
$c: \mathbb R^d \times \mathbb R^d \to \mathbb R^+$ by 
\[
c(x,y)=d_M(b(x),b(y))^2,
\]
where \(d_M\) denotes the Riemannian distance defined in (\ref{dm}). 
Assume that for some $x_0\in \mathbb R^d$, 
\[
\int_{\mathbb{R}^d} d_M(b(x),b(x_0))^2\,d\mu_0(x)<\infty 
\quad \text{and}\quad 
\int_{\mathbb{R}^d} d_M(b(y),b(x_0))^2\,d\mu_1(y)<\infty, 
\]
and that  $\mu_0$  is absolutely continuous with respect to the Lebesgue measure on $\mathbb{R}^d$. 

Then there exists a unique optimal transport map
$
r^\star:\mathbb{R}^d\rightarrow\mathbb{R}^d
$
pushing \(\mu_0\) forward to \(\mu_1\). Moreover, the induced map 
\[
T^\star=b\circ r^\star\circ b^{-1}:  b(\mathbb R^d)   \rightarrow M,
\]
is characterized by
\[
T^\star(\xi)
=
\exp_{\xi}\bigl(-\nabla_M\phi^\star(\xi)\bigr),
\]
where \(\phi^\star\) is a \(c\)-convex Kantorovich potential on \(M\), and $\exp_\xi$  is the Riemannian exponential map at  $\xi \in M$. Equivalently, the optimal transport map on $\mathbb R^d$  is given explicitly by
\[
r^\star(x)
=
b^{-1}\left(
\exp_{b(x)}
\bigl(-\nabla_M\phi^\star(b(x))\bigr)
\right).
\]
\end{thm}

\begin{proof}
Let
\[
\nu_0=b_{\#}\mu_0,\qquad \nu_1=b_{\#}\mu_1
\]
be the push-forward measures on \(M\). Since \(b\) is injective, the transport problem between
\(\mu_0\) and \(\mu_1\) with cost
\[
c(x,y)=d_M(b(x),b(y))^2
\]
is equivalent to the Riemannian optimal transport problem between \(\nu_0\) and \(\nu_1\) with cost
\[
\tilde c(\xi,\eta)=d_M(\xi,\eta)^2 .
\]
Indeed, for any transport plan
$\pi\in\Pi(\mu_0,\mu_1)$, 
the push-forward
$\gamma=(b,b)_{\#}\pi$
belongs to \(\Pi(\nu_0,\nu_1)\). Conversely any
\(\gamma\in\Pi(\nu_0,\nu_1)\) can be pulled back through \(b^{-1}\). Moreover,
\[
\int_{\mathbb{R}^d\times\mathbb{R}^d}
d_M(b(x),b(y))^2\,d\pi(x,y)
=
\int_{M\times M}
d_M(\xi,\eta)^2\,d\gamma(\xi,\eta).
\]
Therefore, the two optimal transport problems have the same value.

Since b is a smooth embedding,
$\nu_0=b_{\#}\mu_0 \ll \mathrm{vol}_M$.
Hence Fathi-Fagalli's theorem \cite[Thm 1.1, Thm 4.3]{FF10} applies: 
there exists a unique optimal transport map
$T^\star:M\rightarrow M$ 
pushing \(\nu_0\) forward to \(\nu_1\). Moreover, there exists a \(c\)-convex Kantorovich potential
\(\phi^\star\) such that
\[
T^\star(\xi)
=
\exp_\xi\bigl(-\nabla_M\phi^\star(\xi)\bigr)
\]
for \(\nu_0\)-almost every \(\xi\in M\).  We now pull this map back to \(\mathbb{R}^d\). Define
\[
r^\star=b^{-1}\circ T^\star\circ b .
\]
Since \(b\) is injective, \(r^\star\) is well defined. Furthermore,
\[
\begin{aligned}
r^\star_{\#}\mu_0
=
(b^{-1})_{\#}T^\star_{\#}b_{\#}\mu_0
=
(b^{-1})_{\#}T^\star_{\#}\nu_0
=
(b^{-1})_{\#}\nu_1
=
\mu_1 .
\end{aligned}
\]
Hence \(r^\star\) is an admissible transport map from \(\mu_0\) to \(\mu_1\).

Finally, using the characterization of \(T^\star\), we obtain
\[
\begin{aligned}
b(r^\star(x))
&=
T^\star(b(x)) =
\exp_{b(x)}
\bigl(-\nabla_M\phi^\star(b(x))\bigr),
\end{aligned}
\]
or equivalently,
\[
r^\star(x)
=
b^{-1}\left(
\exp_{b(x)}
\bigl(-\nabla_M\phi^\star(b(x))\bigr)
\right).
\]
The uniqueness of \(r^\star\) follows from the uniqueness of \(T^\star\) and the injectivity of
\(b\).
\end{proof}

We next show that some auxiliary results, that are well known for quadratic costs,
remain valid for a cost 
$c(x, y)=d_M(b(x), b(y))^2$ with $d_M$ the geodesic distance on $M$.
We first verify that  $W_A$ defines a metric on the space of probability densities with finite second moments. 
The only non-trivial property to verify is the triangle inequality:  
$$
W_A(\rho_1, \rho_3) \leq  W_A(\rho_1, \rho_2) + W_A(\rho_2, \rho_3) \, , 
 $$
 for any probability densities  $\rho_1, \rho_2, \rho_3$ with finite second moments.
The proof relies on the gluing lemma and the composition of optimal transport maps,  therefore avoiding any smoothness assumptions on the densities  $\{\rho_i\}_{i=1}^3$. 
\begin{lem}
The distance $W_A(\cdot,\cdot)$ satisfies the triangle inequality.
\end{lem}
\begin{proof}
We first consider the case where $\mu_1$ and $\mu_2$ are absolutely continuous, while
$\mu_3$ is arbitrary. Let $r$ be the optimal transport map from $\mu_1$ to $\mu_2$, and
let $s$ be the optimal transport map from $\mu_2$ to $\mu_3$. Since
\[
(s\circ r)_{\#}\mu_1
=
s_{\#}(r_{\#}\mu_1)
=
s_{\#}\mu_2
=
\mu_3,
\]
the composition $s\circ r$ is an admissible transport from $\mu_1$ to $\mu_3$. Therefore,
\[
W_A(\mu_1,\mu_3)^2
\le
\int_{\mathbb R^d}
c\bigl(x,s(r(x))\bigr)\,d\mu_1(x),
\]
where
\[
c(x,y)=d_M\!\bigl(b(x),b(y)\bigr)^2.
\]

Taking square roots and applying the triangle inequality for the Riemannian distance,
\[
d_M\!\bigl(b(x),b(s(r(x)))\bigr)
\le
d_M\!\bigl(b(x),b(r(x))\bigr)
+
d_M\!\bigl(b(r(x)),b(s(r(x)))\bigr),
\]
followed by Minkowski's inequality in $L^2(\mu_1)$, yields
\begin{align*}
W_A(\mu_1,\mu_3)
&\le
\left(
\int_{\mathbb R^d}
d_M\!\bigl(b(x),b(r(x))\bigr)^2\,d\mu_1
\right)^{1/2}
\\
&\qquad
+
\left(
\int_{\mathbb R^d}
d_M\!\bigl(b(r(x)),b(s(r(x)))\bigr)^2\,d\mu_1
\right)^{1/2}.
\end{align*}

Since $r_{\#}\mu_1=\mu_2$, the second term becomes
\[
\left(
\int_{\mathbb R^d}
d_M\!\bigl(b(r(x)),b(s(r(x)))\bigr)^2\,d\mu_1
\right)^{1/2}
=
W_A(\mu_2,\mu_3),
\]
while the first term is exactly $W_A(\mu_1,\mu_2)$. Hence,
\[
W_A(\mu_1,\mu_3)
\le
W_A(\mu_1,\mu_2)
+
W_A(\mu_2,\mu_3).
\]

This proves the triangle inequality when the intermediate measures admit optimal transport
maps. For general probability measures, one first establishes the result for the regularized
measures $\mu_i*\chi_\varepsilon$, which are absolutely continuous, and then passes to the
limit as $\varepsilon\to0$ using \cite[Lemma~5.2]{Sa15}. Therefore, $W_A$ satisfies the
triangle inequality in full generality.
\end{proof}

\begin{lem}\label{lem3.5}
Let $w$ be a probability density on $\mathbb{R}^d$ satisfying
\[
\int_{\mathbb{R}^d} d_M(b(x),b(x_0))^2 w(x)\,dx<\infty
\]
for some fixed $x_0\in\mathbb{R}^d$. Then the functional
\[
\rho\longmapsto W_A(\rho,w)^2
\]
is strictly convex on the space of probability measures with finite second
moment with respect to the induced Riemannian distance.
\end{lem}

\begin{proof}
We argue by contradiction. Suppose that there exist two distinct probability
measures $\rho_0\neq\rho_1$ and $s\in(0,1)$ such that
\[
W_A(\rho(s),w)^2
=
(1-s)W_A(\rho_0,w)^2
+
sW_A(\rho_1,w)^2,
\]
where
\[
\rho(s)=(1-s)\rho_0+s\rho_1 .
\]

Let $\gamma_0$ and $\gamma_1$ be optimal transport plans from $\rho_0$ and
$\rho_1$ to $w$, respectively. Since $w$ is absolutely continuous, the
optimal transport map from $w$ to $\rho_i$ is unique by the Riemannian
optimal transport theorem of Fathi-Fagalli \cite[Theorem 1.1]{FF10}.   Equivalently, after
reversing the direction of transport, we may write
\[
\gamma_i=(T_i,\mathrm{id})_{\#}w,
\qquad i=0,1,
\]
where $T_i$ is the unique optimal transport map pushing $w$ forward to
$\rho_i$.

Define the convex combination of transport plans
\[
\gamma(s)=(1-s)\gamma_0+s\gamma_1 .
\]
Then
\[
\gamma(s)\in\Pi(\rho(s),w),
\]
and hence it is an admissible transport plan between $\rho(s)$ and $w$.
Therefore,
\begin{align*}
W_A(\rho(s),w)^2
&\leq
\int_{\mathbb R^d\times\mathbb R^d}
d_M(b(x),b(y))^2\,d\gamma(s)(x,y)
\\
&=
(1-s)
\int_{\mathbb R^d\times\mathbb R^d}
d_M(b(x),b(y))^2\,d\gamma_0(x,y)
\\
&\quad
+s
\int_{\mathbb R^d\times\mathbb R^d}
d_M(b(x),b(y))^2\,d\gamma_1(x,y)
\\
&=
(1-s)W_A(\rho_0,w)^2
+sW_A(\rho_1,w)^2 .
\end{align*}
By the assumed equality, all inequalities above are equalities. Hence
$\gamma(s)$ is itself an optimal transport plan between $\rho(s)$ and $w$.

However, by the uniqueness of optimal transport maps from the absolutely
continuous measure $w$, every optimal transport plan from $w$ to $\rho(s)$
must be concentrated on the graph of a single map. On the other hand,
\[
\gamma(s)
=
(1-s)(T_0,\mathrm{id})_{\#}w
+
s(T_1,\mathrm{id})_{\#}w
\]
is a convex combination of two graph measures. Such a measure is induced by
a transport map only if
\[
T_0=T_1
\qquad w\text{-a.e.}
\]
Consequently,
\[
\rho_0=(T_0)_{\#}w=(T_1)_{\#}w=\rho_1,
\]
which contradicts the assumption $\rho_0\neq\rho_1$.
Therefore, $\rho\longmapsto W_A(\rho,w)^2$ is strictly convex.
\end{proof}

 Define the admissible class of probability densities by
\begin{equation}\label{admclass}
\mathcal K 
=
\left\{
\rho:\mathbb{R}^d\to[0,\infty)\ \text{measurable}:
\int_{\mathbb{R}^d}\rho(x)\,dx=1,\quad
Q_M(\rho)<\infty
\right\},
\end{equation}
where for a fixed reference point $x_0 \in \mathbb{R}^d$, 
\begin{equation}\label{Qm}
Q_M(\rho) =\int_{\mathbb{R}^d}
d_M(b(x),b(x_0))^2\rho(x)\,dx.
\end{equation}

\begin{lem}\label{lem3.6}
Suppose that $\{\rho^n\}\subset \mathcal K  $ converges weakly in
$L^1(\mathbb{R}^d)$ to $\rho$, and that
\[
\sup_{n\geq 1}  Q_M(\rho^n)<\infty .
\]
Then $\rho\in \mathcal K $. Moreover,   for every $\rho_0\in  \mathcal K $,
\begin{equation}\label{lsc} 
\liminf_{n\to\infty}
W_A(\rho_0,\rho^n)^2
\geq
W_A(\rho_0,\rho)^2 .
\end{equation}
\end{lem}

\begin{proof}
The finiteness of the intrinsic second moment follows from the weak lower
semicontinuity of the weighted $L^1$ norm:
\[
Q_M(\rho)
\leq
\liminf_{n\to\infty}Q_M(\rho^n)
<\infty .
\]
Hence $\rho\in K$. For each $n$, let $\pi^n\in\Pi(\rho_0,\rho^n)$ be an optimal transport plan
for the cost
\[
c(x,y)=d_M(b(x),b(y))^2 .
\]
Since
\begin{align}
\int_{\mathbb R^d}
d_M(b(x),b(x_0))^2\,d\rho_0(x) &< \infty \, ,
\\
\sup_n
\int_{\mathbb R^d}
d_M(b(y),b(x_0))^2\,d\rho^n(y) &< \infty \, .  
\label{measure-bound}
\end{align}
By  \eqref{measure-bound} and the compactness of the sublevel sets of $y\to d_M(b(y), b(x_0))$, the family $\{\rho^n\}$ is tight. 
Since the first marginals of  $\pi_n$  are all equal to $\rho_0$ and the second marginals form a tight family, the
sequence $\{\pi^n\}$ is tight in $\mathbb R^d \times \mathbb R^d$. 
Therefore, up to a subsequence,
\[
\pi^n\rightharpoonup\pi
\]
narrowly for some probability measure $\pi$.  By the convergence of the marginals,
\[
\pi\in\Pi(\rho_0,\rho).
\]
Since the Riemannian distance is continuous, the cost
\[
c(x,y)=d_M(b(x),b(y))^2
\]
is continuous and bounded from below. Hence, by the lower semicontinuity of
the Kantorovich functional,
\[
\liminf_{n\to\infty}
\int_{\mathbb{R}^d\times\mathbb{R}^d}
d_M(b(x),b(y))^2\,d\pi^n(x,y)
\geq
\int_{\mathbb{R}^d\times\mathbb{R}^d}
d_M(b(x),b(y))^2\,d\pi(x,y).
\]

Because $\pi^n$ is optimal,
\[
W_A(\rho_0,\rho^n)^2
=
\int_{\mathbb{R}^d\times\mathbb{R}^d}
d_M(b(x),b(y))^2\,d\pi^n(x,y),
\]
and therefore,
\begin{align*}
\liminf_{n\to\infty}
W_A(\rho_0,\rho^n)^2
&\geq
\int_{\mathbb{R}^d\times\mathbb{R}^d}
d_M(b(x),b(y))^2\,d\pi(x,y)
\\
&\geq
W_A(\rho_0,\rho)^2 .
\end{align*}
This completes the proof.
\end{proof}



\section{Variational time-discretization}\label{sec:variation}\label{Sec4} 

We define a time stepping scheme, by first taking a small step size $\tau>0$ and an
initial profile $\rho^0= \rho_0 \in \mathcal K $, 
 where  $\mathcal K $ is the admissible class of densities  defined in 
 \eqref{admclass}.
We  inductively  define $\{\rho^n\}_{n=1}^\infty \subset \mathcal K $, by selecting  $\rho^{n+1}\in \mathcal K $  to be the minimizer of the functional
\begin{align}\notag
\rho^{n+1} = \underset{\rho \in K}{{\rm argmin}} \left ( \frac{1}{2\tau} W_A(\rho^{n}, \rho)^2 +F(\rho) \right ) \, .
\end{align}
The updated density $\rho^{n+1}$  strikes a balance between minimizing the free energy $F(\rho)$ and  the distance $\frac{1}{2\tau} W_A(\rho, \rho^{n})^2$ 
to the density at step $n$. 

This iterative scheme is called the JKO scheme and was originally set forth in \cite{JKO98} for the linear Fokker-Planck equation with mobility $A=I$.  The modifications that are necessary here relate to the use of the generalized  Wassertsein metric $W_A(\rho, \rho^n)$. 
\subsection{The a-priori estimates} 
We set 
$$
G_n (\rho) := \frac{1}{2\tau} W_A(\rho^{n}, \rho)^2 +F(\rho)
$$ 
and suppose $F(\rho^0)<\infty$. 
Since $\rho^{n+1}$ minimizes $G_n(\rho)$, we have 
$$
F(\rho^{n+1})+\frac{1}{2\tau}W_A(\rho^n, \rho^{n+1})^2\leq F(\rho^n).
$$
As a consequence we obtain an energy estimate,
\begin{align}\label{en}
 \sup_{n\in \mathbb{N}} F(\rho^n)\leq F(\rho^0), 
\end{align}
 and 
\begin{align}\label{sq}
\frac{1}{2\tau}\sum_{n\in \mathbb{N}} W_A(\rho^n, \rho^{n+1})^2 \leq F(\rho^0)- \inf F.
\end{align}
Note that this estimate is useful only if $F(\rho_0)<\infty$  and $F$ is bounded from below. The former condition is a natural restriction to finite-energy initial data, and the latter is a reasonable assumption which holds true for many choices of the energy $F(\rho)$.  

From the total square estimate (\ref{sq}) we deduce the usual $1/2$-H\"{o}lder-estimate in time. 
Indeed, for any $0 \leq m \leq n$, the numerical solution $\rho^n$ satisfies 
\begin{align}
W_A(\rho^m,\rho^n)
&\leq \sum_{k=m}^{n-1}
W_A(\rho^k,\rho^{k+1}) \notag\\
&\leq \left[
(n-m)\sum_{k=m}^{n-1}
W_A^2(\rho^k,\rho^{k+1})
\right]^{1/2} \notag\\
&\leq \sqrt{2\tau(n-m)}
\left[F(\rho^0)-\inf F\right]^{1/2}.
\label{mn}
\end{align}
\red{The variational time discretization can be formulated for general energy functionals $F$. 
However, the convergence analysis is restricted to a specific class of energy functionals and requires suitable assumptions on the initial data.}
From now on we restrict to energies
$$
F(\rho)=\int U(\rho, x)dx
$$
with 
$$
U(z, x)= h (z) +z \Psi,
$$
where $h : \mathbb{R}^+ \to \mathbb{R}$ is a density of internal energy, and $\Psi: \mathbb{R}^d \to \mathbb{R}$ is a confinement potential function.  
For $h (z)=z \log z$ and a class of $\Psi \in C^\infty(\mathbb{R}^d)$ satisfying 
$$
\Psi \geq 0, \quad |\nabla \Psi|\leq C(|\Psi|+1) \; \text{for all}\; x\in \mathbb{R}^d,
$$
and $A=I$, the scheme and its convergence theory has been established in \cite{JKO98}.  In the present case of general $A$,  convexity and weak convergence arguments 
still enable us to conclude the following.
\begin{prop}\label{prop1}
Fix $\tau>0$. Let $\rho_0 \in \mathcal K  $ satisfy $F(\rho_0)<\infty$, then there exists a unique $\rho^{*}\in K$ minimizing the functional
$G_n(\rho)$ over $\mathcal K $. That is, for $\rho^{n+1}=\rho^{*}$, 
\begin{equation}\label{gg}
\rho^{n+1}= {\rm argmin}_{\rho \in K}G_n(\rho), \; n=1, 2, \dots  
\end{equation}
\end{prop}

\smallskip
\noindent
The proof of existence for the minimizers is deferred to Section \ref{sec5}. We establish next the Euler-Lagrange equations.

\subsection{Euler-Lagrange equation}
Let
$r(x): \mathbb{R}^d \to \mathbb{R}^d$ be an invertible map defined $\rho^n-a.e.$ with $r_{\#}(\rho^n \mathcal{L}^d)=\rho\mathcal{L}^d$ in the sense that 
$$
\int_{\mathbb{R}^d} f(r(x))\rho^n(x)dx=\int_{\mathbb{R}^d} f(x)\rho(x)dx
$$
for all $f\in C_b({\mathbb{R}^d})$. We then have 
$$
\rho(r(x))=\frac{\rho^n(x)}{\det(\nabla_x r(x))},
$$
where $\nabla_x r$ is the Jacobian matrix of the mapping $r(x)$, and using this
\begin{align*}
F(\rho) & =\int_{\mathbb{R}^d} U(\rho(y), y)dy \\
           & =\int_{\mathbb{R}^d} U\left(\frac{\rho^n(x)}{\det(\nabla_x r(x))}, r(x)\right)det( \nabla_x  r(x))dx,
\end{align*}
which is now a functional of $r$. 

Let $r^{n+1}(x)$ be an optimal map so that 
\begin{equation}\label{risoptimal}
r^{n+1}_{\#}(\rho^n \mathcal{L}^d)=\rho^{n+1}\mathcal{L}^d,
\end{equation}
then we have the following.

\begin{lem} For $\xi \in C_0^\infty(\mathbb{R}^d; \mathbb{R}^d)$ a smooth, compactly supported vector field, we have
\begin{align} \notag
& \frac{1}{2\tau} \int_{\mathbb{R}^d} 
\partial_y c((r^{n+1})^{-1}, x)\cdot \xi(x)\rho^{n+1}(x)dx \\ \label{ws}
& +
\int_{\mathbb{R}^d} \left( -\alpha(\rho^{n+1}, x) {\rm div}(\xi)(x)   +\partial_x U(\rho^{n+1}(x), x))\xi(x)\right)dx=0,
\end{align}
where $\alpha(z, x)=U_z(z, x)z-U(z, x).$ 
\end{lem}
\begin{proof}For any $\xi \in C_0^\infty(\mathbb{R}^d; \mathbb{R}^d)$ and $\epsilon \geq 0$, we define 
$$
r^\epsilon(x):= (id +\epsilon \xi)\circ r^{n+1}(x), 
$$
\red{and let $\rho^{\epsilon}$ be the pushforward of $\rho^n$ by $ r^\epsilon$}, i.e.,  $ r^{\epsilon}_{\#}(\rho^n \mathcal{L}^d)=\rho^{\epsilon}\mathcal{L}^d$. For $\epsilon >0$, $r^\epsilon$ is not 
necessarily optimal, but for any $\epsilon \geq 0$,  
$$
W_A(\rho^n, \rho^\epsilon)^2 \leq \int_{\mathbb{R}^d} c(x, r^\epsilon(x)) \rho^n dx.
$$
Thus 
\begin{align*}
& \limsup_{\epsilon \to 0}\frac{1}{\epsilon}\left[ W_A(\rho^n, \rho^\epsilon)^2-W_A(\rho^n, \rho^{n+1})^2\right]\\
& \leq \limsup_{\epsilon \to 0}\frac{1}{\epsilon}\left[ 
\int_{\red {\mathbb{R}^d}}  \left(    c(x, r^\epsilon(x))  - c(x, r^{n+1}(x))\right)\rho^n dx\right]\\
& =\int_{\mathbb{R}^d}  \xi(r^{n+1})\cdot \partial_y c(x,  r^{n+1}(x)) \rho^n(x)dx \\
& = \int_{\mathbb{R}^d}  \xi(y)\cdot \partial_y c((r^{n+1})^{-1}(y), y)\rho^{n+1}(y)dy.
\end{align*}
For free energy $F(\rho)=\int_{\mathbb{R}^d} U(\rho, x)dx$, we use   
$$
\frac{d}{d\epsilon} \det(\nabla_x  r^\epsilon)\Big|_{\epsilon=0}={\rm div}(\xi) \circ r^{n+1}\det( \nabla_x  r^{n+1}), 
$$
in order to compute   
\begin{align*}
& \lim_{\epsilon \to 0}\frac{1}{\epsilon}
\left[ F(\rho^\epsilon)-F(\rho^{n+1})\right]\\
& = \lim_{\epsilon \to 0}\frac{1}{\epsilon}
 \int_{\mathbb{R}^d} \left[U\left(\frac{\rho^n(x)}{\det(\nabla_x  r^\epsilon(x))}, r^\epsilon(x)\right)\det( \nabla_x  r^\epsilon(x)) \right.\\
& \qquad \left. -U\left(\frac{\rho^n(x)}{\det(\nabla_x  r^{n+1}(x))}, r^{n+1}(x)\right)\det( \nabla_x  r^{n+1}(x))\right]dx 
\\
& = \int_{\mathbb{R}^d} \left[ -\alpha(a(x), r^{n+1}) {\rm div}(\xi) \circ r^{n+1}   +\partial_x U(a(x), r^{n+1})\xi \circ r^{n+1}\right] \det( \nabla_x  r^{n+1})dx,\\
& =\int_{\mathbb{R}^d} \left( -\alpha(\rho^{n+1}, x) {\rm div}(\xi)(x)   +\partial_x U(\rho^{n+1}(x), x)\xi(x)\right)dx,
\end{align*}
\red{ where  $$
a(x): =\frac{\rho^n(x)}{\det(\nabla_x  r^{n+1}(x))} \overset{\eqref{risoptimal}}{=} \rho^{n+1} (r^{n+1} (x)).
$$
}
Hence we have 
\begin{align*}
& 0 \leq \frac{1}{2\tau}  \int_{\mathbb{R}^d} \xi(x)\cdot \nabla_y c((r^{n+1})^{-1}, x)\rho^{n+1}(x)dx  \\
& \qquad +\int_{\mathbb{R}^d} \left( -\alpha(\rho^{n+1}, x) {\rm div}(\xi)(x)   +\partial_x U(\rho^{n+1}(x), x)\xi(x)\right)dx.
\end{align*}
By changing the sign of $\xi$ we obtain the converse inequality, so we actually have equality here.
\end{proof}

The strong form of (\ref{ws}) yields 
\begin{align*}
\frac{1}{2\tau}  \partial_y c((r^{n+1})^{-1}, x) \rho^{n+1}(x)
& =  - \nabla_x \alpha (\rho^{n+1}(x), x) - \partial_x  U(\rho^{n+1}(x), x) \\ 
& = - \rho^{n+1}(x) \nabla_x U_z  (\rho^{n+1}(x), x). 
\end{align*} 

\begin{rem} 
 For the heat equation $\partial_t \rho =\Delta \rho$, we have $b(x)=x$, $c(x, y)=|x-y|^2$,  and $U(z, x)=z {\rm log} z$, the corresponding scheme becomes  
\begin{align*}
& \frac{y-x}{\tau}=-\nabla_y ({\rm log} \rho^{n+1}(y)), \quad \rho^{n+1}(y)=\frac{\rho^n(x)}{\det(\nabla_x  y)}.
\end{align*}
\end{rem}

\blue{ 
We now derive the following estimate for  the Riemannian distance cost.
\begin{lem}\label{dm+} 
Let $\xi=A\nabla_x \phi$,  and define $\psi=\phi\circ b^{-1}$. Then   
$$
\frac{1}{2\tau}
\int_{\mathbb R^d}
\partial_y c((r^{n+1})^{-1}(x),x)
\cdot \xi(x)
\rho^{n+1}(x)\,dx =\frac{1}{\tau}
\int_{\mathbb R^d}
(\rho^{n+1}-\rho^n)\psi\,dx
+
R_n[\psi],  
$$ 
where 
$$
    |R_n[\psi]|
    \leq
    \frac{1}{2\tau}
    \|\nabla_M^2\psi\|_{L^\infty}
    W^2_A(\rho^n, \rho^{n+1}). 
$$
\end{lem} 
} 


\begin{proof}    We use the standard first-variation formula and a second-order Taylor estimate along minimizing Riemannian geodesics. We briefly recall the relevant notation.

For $p\in M$, the exponential map $\exp_p\to M$ maps an initial velocity vector to the endpoint of the corresponding geodesic at time one. Whenever $q$ lies in a normal neighborhood of $p$, we define
$$
\operatorname{Log}_p(q):=\exp_p^{-1}(q)\in T_p M.
$$
Thus, $\operatorname{Log}_p(q)$ is the initial velocity of the minimizing geodesic from $p$ to $q$, parametrized on $[0,1]$.  We write $\nabla_M\psi$ for the Riemannian gradient and $\nabla_M^2\psi$ for the Riemannian Hessian of $\psi$.

\paragraph{\bf Step 1. First variation of the squared distance.}
Away from the cut locus,  the first-variational formula (e.g., \cite[Chapter III, Sec 2, Proposition 2.2]{Sa1996})
gives  
\begin{equation}
    \nabla_q \frac{1}{2}d_M(q,p)^2
    =
    -\operatorname{Log}_{q}(p).
\end{equation}
To derive this gradient identity, let \(\gamma:[0,1]\to M\) be the minimizing geodesic from \(q\) to \(p\), parametrized at constant speed. Varying the initial endpoint \(q\) in the direction \(v\in T_qM\), while keeping \(p\) fixed, and differentiating the geodesic energy gives
\[
\nabla_q\left(\frac12d_M(q,p)^2\right)[v]
=-\langle v,\dot\gamma(0)\rangle_g.
\]
Since \(\dot\gamma(0)=\operatorname{Log}_q(p)\), the desired identity follows.

Consequently, for $
    c(x,y)=d_M(b(x),b(y))^2$,
the chain rule  gives
\begin{equation}\label{cy} 
    \partial_y c(x,y)
    =
    -2(\nabla b(y))^{\top}
    \operatorname{Log}_{b(y)}(b(x)).
\end{equation}
Here $\partial_y$ denotes the Euclidean gradient with respect to $y$, whereas the gradient of the squared distance on $M$ is the Riemannian gradient. \\ 

\paragraph{\bf Step 2. Evaluation of the transport term.}
Substituting (\ref{cy}) 
into the left-hand side of the lemma gives 
$$
\begin{aligned}
  &   \frac{1}{2\tau}
    \int_{\mathbb R^d}
    \partial_y c((r^{n+1})^{-1}(x),x)
    \cdot \xi(x)\rho^{n+1}(x)\,dx  \\ 
& = -\frac{1}{\tau}
\int_{\mathbb R^d}
(\nabla b(x))^{\top}
\operatorname{Log}_{b(x)}
\big(
b((r^{n+1})^{-1}(x))
\big)
\cdot \xi(x)
\rho^{n+1}(x)\,dx .
\end{aligned}
$$
\paragraph{\bf Step 3. Representation via a Riemannian potential.}
Writing $\psi\circ b=\varphi$ and using the chain rule,  the Euclidean gradient of the pullback $\phi(x)=\psi(b(x))$ relates to the ambient gradient via:
$$
\nabla_x \phi(x) =(\nabla b)^\top \nabla_y \psi(b(x)). 
$$
Since $M = b(\R^d)$ is an embedded manifold to $\R^q$, the Riemannian gradient $\nabla_M\psi(b(x))$  is in fact a surface gradient, obtained by projecting the ambient gradient onto the manifold's tangent space.
 The tangent space  $T_p M$ of $M = b(\R^d)$ has as basis the vectors $\del_{x_i} b$, $i = 1, \dots d$. Using this 
 we may compute the orthogonal projector and show that 
$P=\nabla b(x) A(x) (\nabla b)^\top$, where  
$A= ((\nabla b)^\top \nabla b)^{-1}$.
Hence 
\begin{align*}
\nabla_M\psi(b(x)) & =P \nabla_y \psi(b(x)) \\ 
& =\nabla b(x)A(x) ((\nabla b)^\top \nabla_y \psi(b(x)))\\
& =\nabla b(x) A(x) \nabla_x \phi(x) =\nabla b(x)\xi(x)
\end{align*} 
where $\xi=A\nabla_x\phi$. 
Therefore,
\begin{equation}\label{ll} 
\begin{aligned}
&\frac{1}{2\tau}
\int_{\mathbb R^d}
\partial_y c((r^{n+1})^{-1}(x),x)
\cdot \xi(x)
\rho^{n+1}(x)\,dx
\\
&=
-\frac{1}{\tau}
\int_{\mathbb R^d}
\left\langle
\nabla_M\psi(b(x)),
\operatorname{Log}_{b(x)}
\big(
b((r^{n+1})^{-1}(x))
\big)
\right\rangle_g 
\rho^{n+1}(x)\,dx .
\end{aligned}
\end{equation}
\paragraph{\bf Step 4. Taylor expansion along a minimizing geodesic.}
Set
\[
    p=b(r^{n+1}(a)),
    \qquad
    q=b(a).
\]
The second-order Taylor formula along a minimizing geodesic from $p$ to $q$ gives 
\begin{equation}\label{tay2} 
    \psi(q) =
    \psi(p)
    +
    \left\langle
    \nabla_M\psi(p),
    \operatorname{Log}_{p}(q)
    \right\rangle_g
    +
    R(p,q),
\end{equation}
where the remainder satisfies
\begin{equation}\label{rpq} 
    |R(p,q)|
    \leq
    \frac{1}{2}
    \|\nabla_M^2\psi\|_{L^\infty(M)}
    d_M(p,q)^2 .
\end{equation}
The estimate follows by integrating the second derivative of $\psi$ along the geodesic. 
Changing variables \(x=r^{n+1}(a)\) in (\ref{ll}) and applying (\ref{tay2}), we obtain
\begin{equation*}
\begin{aligned}
&\frac1{2\tau}\int_{\mathbb R^d}
\partial_y c(r_{n+1}^{-1}(x),x)
\cdot\xi(x)\rho_{n+1}(x) dx \\ 
&=\frac1{\tau}\int_{\mathbb R^d}
\bigl[\psi(b(r_{n+1}(a)))-\psi(b(a))\bigr]
\rho_n(a) da + R_n[\psi] \\
&=\frac1{\tau}\int_{\mathbb R^d}
(\rho_{n+1}-\rho_n)\psi dx+R_n[\psi],
\end{aligned}
\end{equation*}
where
\begin{equation}
    R_n[\psi]
    =
    \frac{1}{\tau}
    \int_{\mathbb R^d}
    R\big(
        b(r^{n+1}(a)),
        b(a)
    \big)
    \rho^n(a)\,da .
\end{equation}
Finally, (\ref{rpq}) yields 
\begin{equation}
    |R_n[\psi]|
    \leq
    \frac{1}{2\tau}
    \|\nabla_M^2\psi\|_{L^\infty}
    \int_{\mathbb R^d}
    d_M
    \big(
        b(r^{n+1}(a)),
        b(a)
    \big)^2
    \rho^n(a)\,da .
\end{equation}
Since $r_{n+1}$ is optimal for the transport cost induced by $d_M(b(x),b(y))^2$, the integral is $W_{2,A}^2(\rho_n,\rho_{n+1})$. This proves the desired estimate. 
\end{proof}

\begin{lem} \label{lem2} Let $\rho^{n+1}$ be the solution of the iterative scheme (\ref{gg}), then  
$$
\frac{\rho^{n+1}-\rho^n}{\tau}=\nabla_x \cdot ( \rho^{n+1} A(x) \nabla_x U_z(\rho^{n+1}, x)) +T_1, 
$$
\tcc{where the truncation error  
$$
\left|\int T_1 \psi(b(x))dx\right| \leq \frac{1}{2\tau}     \|\nabla_M^2\psi\|_{L^\infty} W_A(\rho^n, \rho^{n+1})^2, 
$$
where $\nabla_M \psi$ and $\nabla_M^2 \psi$ denote the Riemannian gradient and Hessian of $\psi:M\to\mathbb{R}$, respectively, and the pullback function $\psi\circ b$ belongs to $C_0^\infty(\mathbb{R}^d)$.}
\end{lem}
\begin{proof}  Let us fix $\phi \in C_0^\infty (\mathbb{R}^d; \mathbb{R})$ and take $\xi:=A(x) \nabla_x\phi$  in (\ref{ws}). Then the second term is 
\begin{align*}
& II = \int_{\mathbb{R}^d} \left( -\alpha(\rho^{n+1}, x) {\rm div}(\xi)(x)   +\partial_x U(\rho^{n+1}(x), x)\xi(x)\right)dx \\
& = \int_{\mathbb{R}^d} \xi \cdot (\nabla_x \alpha (\rho^{n+1}, x) +\partial_x U(\rho^{n+1}, x))dx \\
& = \int_{\mathbb{R}^d} A\nabla_x \phi \cdot \nabla_x (U_z(\rho^{n+1}, x)) \rho^{n+1} dx\\
& =-\int_{\mathbb{R}^d} \phi(x) \nabla_x \cdot(A \rho^{n+1}\nabla_x U_z(\rho^{n+1}, x))dx.
\end{align*}
We next examine the first term 
$$
I= \frac{1}{2\tau} \int_{\mathbb{R}^d} 
\partial_y c((r^{n+1})^{-1}, x)\cdot \xi(x)\rho^{n+1}(x)dx. 
$$
\tcc{Define $\psi: M\to \mathbb{R}$   such that $\psi=\phi \circ b^{-1}$. Then by  Lemma \ref{dm+}, we obtain 
$$
I= \frac{1}{\tau}
\int_{\mathbb R^d}
(\rho^{n+1}-\rho^n)\psi\,dx
+
R_n[\psi],  
$$ 
where 
$$
    |R_n[\psi]|
    \leq 
    \frac{1}{2\tau}
    \|\nabla_M^2\psi\|_{L^\infty}
    W^2_A(\rho^n, \rho^{n+1}), 
$$} 
which completes the proof.
\end{proof}

\begin{rem} 
The variational scheme is an approximation of the implicit Euler scheme
$$
\frac{\rho^{n+1}(x)-\rho^n(x)}{\tau}=\nabla_x \cdot ( \rho^{n+1} A(x) \nabla_x U_z(\rho^{n+1}, x)),
$$
the approximate error being controlled quadratically in the generalized Wasserstein distance.

\end{rem}

\subsection{Convergence} 
Next, \red{we introduce a piecewise-constant interpolation  $\rho_\tau: [0,\infty) \to K$  defined by setting $t_n= n\tau$ and 
 $$
 \rho_\tau(t) = \rho^n \quad \text{on each time interval} \quad  t \in [t_n, t_{n+1})   \; (n = 0, 1, \cdots).
 $$  } 
Let us suppose now that as $\tau \to 0$,
\begin{equation}\label{limit}
{\rho_\tau \to \rho \; \text{strongly in}\;  L^1_{\rm loc}(\mathbb{R}^d \times (0,\infty)).}
\end{equation}
For the sake of generality,  \red{we assume strong convergence $\rho_\tau \to \rho $ (for instance, pointwise almost everywhere),  ensuring that the nonlinear terms pass to the limit. This is a stringent assumption that must be verified for any specific choice of the functional $F(\rho)$.}   We shall discuss in Section \ref{sec5} a strategy  to retrieve compactness for linear diffusion, for which convergence in $C([0,T],L^1_{\rm weak}(\mathbb{R}^d))$ will suffice.  
\begin{thm}The function $\rho$ is a weak solution of the nonlinear  problem
\begin{align}\label{4.7}
& \partial_t \rho =\nabla_x \cdot (\rho A(x)\nabla_x U_z(\rho, x))\; \text{in}\; \mathbb{R}^d \times (0,\infty),\\
& \rho(x, 0)=\rho_0.
\end{align}
\end{thm}
\begin{proof}
Let $\phi \in C_0^\infty (\mathbb{R}^d; \mathbb{R})$ be arbitrary and choose  $\psi= \phi\circ b^{-1}$. Then by Lemma  \ref{lem2}, we obtain 
\begin{align*}
&  \int_{\mathbb{R}^d} \frac{\rho^{n+1}(x) -\rho^n(x)}{\tau}\phi(x)dx  
+ \int_{\mathbb{R}^d}  \nabla \phi  \cdot(A \rho^{n+1}\nabla_x U_z(\rho^{n+1}, x))dx   =  - R_n[\psi].
\end{align*}
Passing to the piecewise-constant interpolant $\rho^\tau$, it follows that, for every test function $\phi \in  C_0^\infty (\mathbb{R}^d \times [0, \infty);  \mathbb{R})$, 
\begin{align*}
&  \int_{\mathbb{R}^d} \frac{\rho_\tau(t+\tau) -\rho_\tau(t)}{\tau}\phi(x, t)dx 
+ \int_{\mathbb{R}^d} \nabla_x \phi(x, t) \cdot(A \rho_\tau(t+\tau)\nabla_x U_z(\rho_\tau(t+\tau), x))dx  \\
& = - R_n^\tau [\psi], \quad t\in [t_n, t_{n+1}),
\end{align*}
where $R_n^\tau [\psi]$ denotes the corresponding remainder term.  Hence 
\begin{align*}
& \left| -\int_\tau^\infty \int_{\mathbb{R}^d} \frac{\phi(x, t) -\phi(x, t-\tau)}{\tau}\rho_\tau(t)dxdt  -\frac{1}{\tau}\int_0^\tau\int_{\mathbb{R}^d}
\rho^0(x) \phi(x, t)dx dt \right. 
\\  
& \quad \left. + \int_\tau^\infty  \int_{\mathbb{R}^d} \nabla_x \phi(x, t-\tau) \cdot(A \rho_\tau(t)\nabla_x U_z(\rho_\tau(t), x))dx \right|  \\
&  \leq  \sum_{n=0}^\infty \int_{t_n}^{t_{n+1}} |R_n^\tau [\psi(\cdot, t)]|dt  \\
& \leq   \sum_{n=0}^\infty \int_{t_n}^{t_{n+1}} \frac{1}{2\tau}\sup_{M}|\nabla_M^2 \psi(\cdot, t)| W_A(\rho^n, \rho^{n+1})^2 dt \\
& \leq \frac{1}{2} \sup_{M \times \mathbb{R}^+}|\nabla_M^2 \psi(\cdot, \cdot)| \sum_{n=0}^\infty W_A(\rho^n, \rho^{n+1})^2 
\\
&\overset{\eqref{sq}}{\leq} C \tau \, .
\end{align*}
\red{ Passing to the limit as $\tau \to 0$,} 
\begin{align}\label{wk}
\int_0^\infty \int_{\mathbb{R}^d} -\rho \phi_t  + \nabla_x U_z(\rho, x)\cdot A(x) \rho(x, t)\nabla_x \phi dxdt -\int_{\mathbb{R}^d} \rho_0\phi(\cdot, 0)dx=0,
\end{align}
and $\rho$ is a weak solution of the nonlinear equation (\ref{4.7}). 
\end{proof}

\section{A priori estimates and compactness}\label{sec5} 
\subsection{A priori estimates} 
Whether the variational scheme (\ref{gg}) is well-posed also relies on the structure of $F$. 
The typical form of $F$ is 
$$
F(\rho)=\int_{\mathbb{R}^d} U(\rho)dx +\int_{\mathbb{R}^d} \Psi(x) \rho dx+\frac{1}{2}
\int_{\mathbb{R}^d\times \mathbb{R}^d} \Gamma(x-y)\rho(x)\rho(y)dxdy,
 $$
where $\Gamma: \mathbb{R}^d \to \mathbb{R}$ is an interaction potential. According to the model (\ref{main})  this corresponds to 
the aggregation--diffusion equation 
$$
 \partial_t \rho = \nabla \cdot (\rho A(x)\nabla (U'(\rho) + \Psi(x) +\Gamma*\rho)) \, .
 $$
 To simplify the presentation, we restrict our analysis to the case without interactions ($\Gamma=0$) and consider only 
linear diffusion. The extension to more general energy functionals is left for future investigation.
 
 We first establish two fundamental a priori estimates for the linear diffusion case, which will play a central role in the convergence 
 analysis of the variational scheme.  

\begin{lem} \label{lemc1}
(i) For any $\rho \in K$ and $\epsilon>0$, there exists $C_\epsilon>0$ such that 
\begin{align}\label{sm}
S(\rho):= \int \rho \log \rho dx   \geq  -C_\epsilon  -\epsilon Q_M(\rho), 
\end{align}
where $Q_M(\rho)$ is defined in (\ref{Qm}).  \\
(ii) For any $\rho_1, \rho_2 \in K$, we have 
\begin{align}\label{mw}
Q_M(\rho_2) \leq 2W_A^2(\rho_1, \rho_2)+2Q_M(\rho_1).
\end{align}
\end{lem}
\begin{proof} 
(i) {\bf Entropy estimate.} Let $D(x):=d_M(b(x), b(x_0))$ denotes  the intrinsic distance. By assumptions (\ref{c0}) and \eqref{abb}, 
$$
D^2(x) \geq |b(x)-b(x_0)|^2 \geq c_0|x-x_0|^2.  
$$
Define the reference probability density 
$$
\sigma(x)=e^{-2 D(x)}/Z, \quad Z=\int e^{-2 D(x)}dx <\infty,
$$
where the finiteness of $Z$ follows from the above coercivity estimate.

Since $\int \rho dx = \int \sigma dx =1$,  the nonneativity of relative entropy gives 
$$
\begin{aligned}
0 &\le \int  \rho \log \rho - \sigma \log \sigma - (1+ \log \sigma) (\rho - \sigma) dx  
\\
&= \int \rho \log \frac{\rho}{\sigma}dx  
\\
& \leq S(\rho)+2\int D(x) \rho dx  +\log Z
\\
&  \leq S(\rho) +\epsilon Q_M(\rho)+C_\epsilon, \quad C_\epsilon=1/\epsilon +\log Z,
\end{aligned}
$$
where we used the Young's inequality $D(x) \leq \epsilon D(x)^2 +\frac{1}{\epsilon}$. Consequently, 
$$
S(\rho) \geq - \epsilon Q_M(\rho) - C_\epsilon.
$$ 
(ii) {\bf Second-moment estimate.} 
Let $r^*$ be an optimal transport map satisfying $ r^{*}_{\#}(\rho_1 \mathcal{L}^d)=\rho_2 \mathcal{L}^d$.
By the triangle inequality 
 $$
\begin{aligned}
d_M(b(r^*(x)),b(x_0))
&\leq d_M(b(r^*(x)),b(x))\\
&\quad+d_M(b(x),b(x_0)).
\end{aligned}
$$
Hence,
\[
\begin{aligned}
Q_M(\rho_2)
&=\int_{\mathbb R^d}
d_M(b(r^*(x)),b(x_0))^2\rho_1(x)\,dx\\
&\leq 2\int_{\mathbb R^d}
d_M(b(r^*(x)),b(x))^2\rho_1(x)\,dx +2Q_M(\rho_1)\\
&=2W_{A}^2(\rho_1,\rho_2)+2Q_M(\rho_1),
\end{aligned}
\]
where
\[
W_{A}^2(\rho_1,\rho_2)
:=\inf_{\pi\in\Pi(\rho_1,\rho_2)}
\int_{\mathbb R^d\times\mathbb R^d}
d_M(b(x),b(y))^2\,d\pi(x,y).
\]
Thus,
\[
Q_M(\rho_2)\leq
2W_{A}^2(\rho_1,\rho_2)+2Q_M(\rho_1).
\]
\end{proof} 
\noindent{\bf  Proof of Proposition \ref{prop1}}: 
First, we estimate the functional $G_n$ from below. 
By the lower bound for the generalized  Wasserstein distance and the entropy estimate, we obtain
\begin{align*}
G_n(\rho) & \geq \frac{1}{2\tau} W_A^2(\rho^n, \rho) +S(\rho) \\
& \geq 
-\epsilon Q_M(\rho)-C_\epsilon.
\end{align*}
Now, let $\{\rho_k\}$ be a minimizing sequence for $G_n(\rho)$, so that $\{G_n(\rho_k)\}$ is bounded from below and above. 
By the preceding estimates, we obtain 
$$
\int_{\mathbb{R}^d} \rho_k |\log \rho_k|dx  \leq C, \qquad \sup_k Q_M(\rho_k)   \leq C.
$$
The entropy bound prevents concentration. Indeed, for any $\Lambda >1$ 
$$
\int_{\rho_k \geq \Lambda } \rho_k dx \leq \frac{1}{\log \Lambda } \int_{\rho_k \geq \Lambda  } \rho_k \log \rho_k dx \leq \frac{C}{\log \Lambda } 
$$
which tends to zero uniformly in $k$ as $\Lambda  \to \infty$.  

On the other hand, since $D(x) \geq c_0|x-x_0|^2$, we have for any $R>0$   
$$
\int_{|x-x_0| > R} \rho_k (x) dx \le \frac{1}{R^2} \int |x-x_0|^2 \rho_k(x)dx \leq \frac{Q_M(\rho_k)}{c_0 R^2 }\leq 
\frac{C}{c_0 R^2}.
$$
Thus, the sequence $\{\rho_k\}$ is tight. 

Together, the entropy bound and tightness imply uniform integrability on $\mathbb{R}^d$.  
By the Dunford--Pettis theorem, there exists a subsequence,  still denoted by $\{\rho_k\}$, and a probability density $\rho^*$ such that 
$$
\rho_k  \rightharpoonup  \rho^* \qquad \text{weakly in}\quad L^1(\mathbb R^d). 
$$
Moreover, the lower semicontinuity of \(Q_M\) yields
\[
Q_M(\rho^*)\leq \liminf Q_M(\rho_k) \leq C.
\]
Hence, \(\rho^*\in K\).

It remains to show that $\rho^*$  minimizes  $G_n(\rho)$.  The entropy $S$ and the potential energy $\int \Psi(x) \rho dx$  are  weakly  lower semicontinuous in $L^1$ (see \cite{JKO98}). 
while the lower semicontinuity of the generalized Wasserstein distance $W_A^2$ follows from Lemma \ref{lem3.6}.  Consequently,
\[
G_n(\rho^*) \leq \liminf {G_n}(\rho_k) =\inf {G_n}(\rho),
\]
proving the existence of a minimizer. Moreover, by Lemma 3.5 and the convexity of \(F\), the functional \(G_n\) is strictly convex, ensuring the uniqueness of the minimizer at each time step.
\null                   
   \hfill                  
   $\square$


\subsection{Compactness estimates}  
We next derive the compactness estimates needed for the convergence analysis as $\tau \to 0$, as outlined in (\ref{limit}). 
For simplicity, we restrict attention  to  the linear diffusion case. The main ingredients are  entropy bounds ensuring spatial equi-integrability of
 the approximate solutions, and time-regularity estimates derived from the parabolic structure of the variational scheme.

The proof proceeds in two steps. \\
{\bf Step 1. A priori estimates in time.}   
By Lemma \ref{lemc1},   the free energy $F(\cdot)$ is bounded from below. Moreover, the discrete energy dissipation estimate  (\ref{mn}) yields  
$$
 W_A(\rho^m, \rho^{n}) \leq  \sqrt{2(F(\rho^0) - \inf F)}
 \sqrt{(n-m)\tau}.
 $$
For any $0 \leq s \leq t$,  the triangle inequality gives 
\begin{align*} \notag
W_{A}(\rho_\tau(s), \rho_\tau(t)) & \leq  
W_{A}(\rho_\tau(s), \rho^{\lfloor s/\tau\rfloor +1 }) \\
& \qquad +W_{A}(\rho^{\lfloor s/\tau \rfloor +1}, \rho^{\lfloor t/\tau \rfloor+1 })+
W_{A}(\rho^{\lfloor t/\tau\rfloor +1}, \rho_\tau(t)). 
\end{align*}
Since \(\rho_\tau\) denotes the standard piecewise-constant interpolation
\[
\rho_\tau(t)=\rho^{\lceil t/\tau\rceil},
\qquad t>0,
\]
then the first and third terms in the triangle inequality vanish, apart from possible endpoint conventions. In that case, setting
\[
m=\lceil s/\tau\rceil,\qquad n=\lceil t/\tau\rceil,
\]
we obtain directly
\begin{align} \label{12}
 W_{A}(\rho_\tau(s), \rho_\tau(t))  \leq W_A(\rho^m, \rho^{n})  \leq  C(t-s+\tau)^{1/2}. 
\end{align}
This provides the desired approximate \(1/2\)-Hölder continuity in time, uniformly as \(\tau\to0\). \\ 
{\bf Step 2. Compactness.}   The uniform moment and entropy bounds from Lemma \ref{lemc1}
ensure tightness and equi-integrability of
$\{\rho_\tau(t)\}_{\tau>0}$, uniformly for $t\in[0,T]$.
By the Dunford--Pettis theorem, this family is
relatively weakly compact in $L^1(\mathbb R^d)$
at each fixed time.

Moreover, the time-regularity estimate  (\ref{12})  gives
\[
W_A(\rho_\tau(s),\rho_\tau(t))
\le C\sqrt{|t-s|+\tau},
\qquad 0\le s,t\le T.
\]
Combining spatial compactness with this asymptotic
equicontinuity, an Arzel\`a--Ascoli-type compactness
argument yields a subsequence $\tau_j\to0$ and
a limit curve
\[
\rho\in C([0,T];L^1_{\mathrm{weak}}(\mathbb R^d))
\]
such that, for every $\varphi\in L^\infty(\mathbb R^d)$,
\[
\sup_{t\in[0,T]}
\left|
\int_{\mathbb R^d}
(\rho_{\tau_j}(t,x)-\rho(t,x))
\varphi(x)\,dx
\right|\longrightarrow0.
\]
Thus, $\rho_{\tau_j}$ converges uniformly in time
in the weak $L^1$ topology.

\bigskip
\subsection*{Acknowledgement} 
This work resulted from a collaboration initiated during H. Liu's visit to KAUST. \blue{We thank Professor Julio Valencia Guevara (Universidad Nacional de San Agustín, Peru)  for his insightful feedback and  the anonymous referees for their constructive comments,  which helped improve the manuscript.} H. Liu was partially supported by the National Science Foundation under Grant DMS-1812666. A. E. Tzavaras was supported by King Abdullah University of Science and Technology (KAUST), baseline funds No. BAS/1/1652-01-01.


\bibliographystyle{abbrv}

\appendix
\section{Variational Discretizations and the high-friction limit in Lagrangian coordinates}\label{sec:appA}
Consider an evolution equation of the form
\begin{equation}\label{pp0}
 \ddot Y = -  \delta_Y \cE -  \frac{1}{\zeta} B(x) \dot Y. 
\tag{P1}
\end{equation}
In \eqref{pp0},  the evolution of the motion $Y : (0, \infty) \times \R^d \to \R^d$ is
driven by a functional $\cE(y)$  and $\delta_y \cE$  stands for the (generator of) the functional derivative of $\cE$ 
defined by
$$
\frac{d}{d\tau} \cE ( y + \tau \varphi) \Big |_{\tau = 0} =: \langle \delta_y \cE , \varphi \rangle. 
$$
The evolution is subjected to a heterogeneous frictional environment captured by the  
positive definite and symmetric matrix $B(x) > 0$ with frictional parameter $\tfrac{1}{\zeta}$.

In the theory of relaxation approximations  it is classical to turn this problem to a small-mass approximation limit.
For convenience set $\zeta = \sqrt{\eps}$ and introduce the change of time scale
$
y(x, t) = Y \Big (x, \frac{t}{\sqrt{\eps}} \Big )
$
that captures the long-time response of $Y$. The change of time scale transforms \eqref{pp0} to the limiting problem
\begin{equation}\label{pp1}
\eps \ddot y = -  \delta_y \cE - B(x) \dot y \, .
\tag{P1}
\end{equation}
Its formal small-mass approximation limit $\eps \to 0$ reads
\begin{equation}\label{pp2}
\dot y = -  A(x)  \delta_y \cE   \, ,
\tag{P2}
\end{equation}
where $A(x) = B^{-1} (x)$ is the inverse of $B(x)$.
The problem \eqref{pp1} is endowed with the energy dissipation identity
\begin{equation}\label{energy1}
\frac{d}{dt} \left (  \int_{\mathbb{R}^d}   \eps \tfrac{1}{2} | \dot y |^2 dx + \cE (y)  \right ) + \int_{\mathbb{R}^d}   \dot y \cdot B(x) \dot y \, dx = 0
\end{equation}
We assume here that the potential energy $\cE$ is assumed to be convex and coercive so that \eqref{energy1} describes energy dissipation and is maximally stabilizing. (Relaxing the assumption $\cE$ is also a meaningful problem and leads to issues of relaxations of functionals; See \cite{MO08}
and references therein).

The gradient flow \eqref{pp2} is the formal limit of \eqref{pp1};
the matrix $A(x)$ is defined by $A : = B^{-1}$ and is positive definite and symmetric. The gradient flow 
inherits form  \eqref{energy1} the  identity
\begin{equation}\label{energy2}
\frac{d}{dt} \cE (y)  + \int_{\mathbb{R}^d}   \delta_y \cE \cdot A(x) \delta_y \cE \, dx = 0
\end{equation}
capturing the (potential) energy dissipation identity for \eqref{pp2}.

In this appendix, we present the formalism of the time-step discretization method for the problem \eqref{pp1} (see \cite{Dem97, DST00}) 
and compare with the method of minimizing movements of De Giorgi ({\it e.g.} \cite{Gi93}) for the problem \eqref{pp2}, following ideas 
traced to \cite{MO08}.  
Consider the initial value problem \eqref{pp1} with initial data
\begin{equation}\label{indata1}
y(0, x) = y_0 (x) \, , \quad \frac{\del y}{\del t} (0, x) = v_0 (x)  \, ,
\end{equation}
where $(y_0, v_0)$ are given functions. Consider a discretization of the time interval $(0,T)$ into discrete points $t_j = j \tau$, $j = 0, ... , N$
 with $N \tau = T$.
The method of time-step discretization
consists of approximating $y(t,x)$ by the iterates $\{ y_j (x) \}_{j =1 , ..., N}$ obtained by solving the discretized problems
\begin{equation}\label{ELeqn}
\begin{aligned}
\eps \frac{v_j - v_{j-1}}{\tau} &= -  \delta_y \cE (y_j)  - B(x) v_j
\\
v_j &= \frac{y_j - y_{j-1}}{\tau}
\end{aligned}
\qquad j = 1,\cdots, N.  
\end{equation}
The strategy is to solve \eqref{ELeqn} for $j=1, ... , N$ and pass to the limit $\tau \to 0$ to solve \eqref{pp1} with data \eqref{indata1}.
This is based on the observation that \eqref{ELeqn} are the Euler-Lagrange equations for the constrained variational problem
\begin{equation}\label{vp1}
\min_{ \big \{ (y, v) \; :  \; \; v = \frac{y - y_{j-1}}{\tau} \big \}}  \int_{\mathbb{R}^d}   \tfrac{\eps}{2} | v - v_{j-1} |^2 \, dx  + \cE (y) + \int_{\mathbb{R}^d}   \tfrac{\tau}{2} v \cdot B(x) v \, dx \,, j = 1, \cdots, N. 
\end{equation}
Under convexity and coercivity conditions for $\cE(y)$, the variational problem \eqref{vp1} amounts to  minimizing a convex functional 
subject to an affine contraint. 
The minimizer will satisfy the Euler-Lagrange equations, which we now compute: Let $(y_j,v_j)$ be a minimizer of \eqref{vp1} and consider the 
variation  $y_j + \delta Y$, $v_j + \delta V$ of  $(y_j, v_j)$, with $(Y, V)$ smooth. Then  the constraint implies $Y = \tau V$ while
 taking the variation of the functional leads to
$$
\int_{\mathbb{R}^d}   \eps (v-v_{j-1}) \cdot V + \langle \delta_y \cE (y)  , Y \rangle  + \int_{\mathbb{R}^d}   \tau B(x) v \cdot V \, dx = 0 \, ,
$$
that is $(y_j,v_j)$ satisfies \eqref{ELeqn}. The process produces an approximate solution of \eqref{pp1}. Passing to the limit $\tau \to 0$
is in general a difficult problem tied to the existence theory and compactness properties of the limit equation; see \cite{DST00} for an
example.

Next, we note that the limit  $\eps \to 0$ of the variational problem \eqref{vp1}  leads to the minimization problem 
\begin{equation}\label{vp2}
\min_{ y }  \left ( \cE (y) + \int_{\mathbb{R}^d}  \tfrac{1}{2 \tau}  (y-y_{j-1}) \cdot B(x) (y-y_{j-1}) \, dx \right ) \, .
\end{equation}
The Euler-Lagrange equation for \eqref{vp2} is computed as
$$
\int_{\mathbb{R}^d}   B(x) \tfrac{1}{\tau} (y - y_{j-1}) \cdot Y  \, dx + \langle \delta_y \cE (y), Y \rangle = 0. 
$$
That is, the sequence of minimizers $\{ y_j \}_{j = 1, ... , N}$ of the consecutive problems \eqref{vp2}, $j =1, ... , N$, satisfy
\begin{equation}
\frac{y_j - y_{j-1}}{\tau} = - A(x) \delta_y \cE (y_j) 
\end{equation}
and yield the minimizing movement method for  \eqref{pp2} with initial data $y = y_0(x)$.

A special case for \eqref{pp1} arises for  the one-dimensional case $y : (0,T) \times \R \to \R$ and  the energy
$$
\cE (y) = \int_{\mathbb{R}^d}   W(\del_x y ) dx
$$
with $W(u)$ convex and $\sigma (u) = W^\prime (u)$ monotone. In that case \eqref{pp1} takes the form
of the nonlinear wave equation
\begin{equation}
\eps y_{tt} = \del_x \big( \sigma (y_x) \big )  - b(x) y_t
\end{equation}
with $b(x) > 0$ scalar function.
Under periodic boundary conditions it can be proved following \cite{DST00}, using the methodology  of compensated compactness,
that the time step discretizations generated by the scheme \eqref{ELeqn} converge in the limit $\tau \to 0$ and for $\eps$ fixed to an entropy weak
solution of  the system of conservation laws 
$$
\begin{aligned}
u_t &= v_x 
\\
\eps v_t &= \sigma(u)_x - b(x) v
\end{aligned}
$$
where $v = y_t$, $u = y_x$,  which dissipates all the convex entropies.

\section{The isometric embedding theorems}\label{sec:isoemb}
Here we connect the structural condition
\begin{equation}\label{isoeq}
(\nabla_x b)^\top \nabla_x b= B (x)
\end{equation}
to the isometric embedding theorem.  Given a symmetric, positive definite matrix  $A(x)$, let $B(x) = A^{-1} (x)$ be its inverse.
The idea is to connect  \eqref{isoeq} to the isometric embedding problem  for Riemannian manifolds:  finding a smooth map $b: \mathbb{M}^d \to \mathbb{R}^q$ such that 
$$
db\cdot db =g, 
$$
where the metric $g$ is given by 
$$
g=\sum_{i, j=1}^dG_{i, j}(x)dx_idx_j, \quad G=B .
$$
The system contains $s_d=d(d+1)/2$ (called the Janet dimension) equations, so one can not do better than $q = s_d$ (see Han--Hong \cite{HH06} for a thorough
introduction to the problem of isometric embedding). 

Existence of a smooth local isometric embedding into $\mathbb{R}^{s_d+d}$ is guaranteed from \cite{GR70}.
%
%
%
The smooth global isometric embedding is often called the Nash Embedding Theorem \cite{Na56}. The result was improved by G\"{u}nther \cite{Gu91}
for  $q=\max\{s_d+2d, s_d+d+5\}$; here $q$ is smaller than that in \cite{Na56}. The proof in \cite{Gu91} consists of two steps. First, for any Riemannian metric $g$ on the manifold $\mathbb{M}^d$, one finds an embedding $b_0 : \mathbb{M}^d \to \mathbb{R}^q$ such that $g-db_0 \cdot db_0$  is also a metric. Second, one modifies $b_0$ to get $b$ satisfying $db\cdot db = g$. During this process, it is maintained that $b$ is still an embedding. 

In the present setting,  with $ \mathbb{M}^d$ being the space $\mathbb{R}^d$, endowed with the Riemannnian distance induced by $B $,   we state the following 
\begin{prop} Let $B $ be a smooth $d\times d$ matrix,  symmetric and positive definite. 
Then there exists a smooth isometric embedding map 
$b: \mathbb{M}^d \to \mathbb{R}^q$  with $q=\max\{s_d+2d, s_d+d+5\}$ such that 
$$
(\nabla_x b)^\top \nabla_x b = B (x) \qquad \forall  x\in  \mathbb{M}^d. 
$$
If $q\geq d+1$, then $b\in C^{1}(\mathbb{M}^d; \mathbb{R}^q)$.  
\end{prop} 
The smooth isometric embedding with index $q=\max\{s_d+2d, s_d+d+5\}$ is due to G\"{u}nther \cite{Gu91} and the $C^1$ isometric 
embedding with $q\geq d+1$ is implied by the Nash-Kuiper theorem \cite{Na54, Ku55, Ku59}.  
In recent years, a sequence of papers (see De Lellis and Sz\'{e}kelyhidi \cite{DS17} for a survey) have appeared as  an outgrowth of the celebrated results of Nash and Kuiper for the $C^1$-isometric embeddings of a $d$-dimensional Riemannian manifold into $\mathbb{R}^{m}$, $m \ge d+1$.  
These works have sharpened   the regularity of embeddings in the Nash-Kuiper theorem (see \cite{CLS12, DIS18,CHI25}). 
The central technique is the use of convex integration and the $h$-principle, introduced by M. Gromov \cite{Gr17} as a generalization of the Nash--Kuiper paradigm.
In \cite{CLS12}, using the method of convex integration, Conti-De Lellis, Sz\'{e}kelyhidi Jr. provide a constructive  proof of a version of the Nash--Kuiper theorem. The issue of the best H\"{o}lder exponent $\alpha$ is studied in \cite{CHI25}.
For easy reference, here we record the main result in \cite{CLS12}.
\begin{thm} Let $\mathbb{M}^d$ be an d-dimensional compact manifold with a Riemannian 
metric $g$ in $C^\beta(\mathbb{M}^d)$ and $q \geq d+1$. 
Then there is a constant $\delta_0>0$ such that if $u\in C^2(\mathbb{M}^d; \mathbb{R}^q)$ 
and $\alpha$ satisfy
$$
\|\partial_i u\cdot \partial_j u-g_{ij}\|_0\leq \delta_0^2\; \text{and} \; \; 0<\alpha<\min\left\{
\frac{1}{d(d+1)^2}, \frac{\beta}{2} \right\}, 
$$
there exists a map $v$ in $C^{1,\alpha}(\mathbb{M}^d;\mathbb{R}^q)$ such that 
$$
\partial_i v\cdot \partial_j v=g_{ij} \; \;\text{and} \; \;\|v-u\|_1 \leq C \delta_0^2.
$$
\end{thm}

\end{document}